\documentclass[preprint,12pt,sort&compress]{elsarticle}
\usepackage{amsmath}
\usepackage{mathtools}
\usepackage{amssymb}
\usepackage[mathscr]{eucal}
\usepackage{amscd}
\usepackage{color}
\usepackage{verbatim}
\usepackage{sistyle}
\usepackage{subfigure}
\usepackage{makecell}
\usepackage[usenames,dvipsnames]{xcolor}
\usepackage{tikz-cd}
\usepackage{algorithm}
\usepackage{xfrac}
\usepackage[permil]{overpic}
\DeclareMathOperator{\SO3}{{{\rm SO(3)}}} 
\DeclareMathOperator{\so3}{{{\rm so(3)}}}

\newcommand{\f}[1]{{\boldsymbol{#1}}}
\newcommand{\baf}[1]{{\bar{\boldsymbol{#1}}}}

\newcommand{\sub}{\subset}
\newcommand{\bEq}{\begin{equation}}
\newcommand{\eEq}{\end{equation}}
\newcommand{\beq}{\begin{equation*}}
\newcommand{\eeq}{\end{equation*}}

\newcommand{\car}{\times}
\newcommand{\mto}{\mapsto}

\newcommand{\byd}{\,{\raisebox{.092ex}{\rm :}{\rm =}}\,}

\newcommand{\fR}[1]{{\mathbf{#1}}}
\newcommand{\hfR}[1]{\hat{\mathbf{#1}}}

\newcommand{\sepr}[1]{\,\,\,\,\textnormal{{#1}}\,\,\,\,}

\newcommand{\lf}{\left}
\newcommand{\rg}{\right}

\newcommand{\sS}[1]{{\scriptscriptstyle {#1}}}
\newcommand{\Tra}{^{\mathsf{\sS\!T}}}
\newcommand{\Rn}{\text{I\!R}}
\newcommand{\tif}[1]{{\widetilde{\boldsymbol{#1}}}}

\newcommand{\bAl}{\begin{align}}

\newcommand{\Gam}{\varGamma}
\newcommand{\B}[1]{{\mathbb{#1}}}
\newcommand{\veps}{\varepsilon}
\newcommand{\gam}{\gamma}
\newcommand{\del}{\delta}

\newcommand{\vtht}{\vartheta}
\newcommand{\vTht}{\varTheta}

\newcommand{\fr}[2]{\frac{#1}{#2}\,}

\newcommand{\wti}[1]{{\widetilde{#1}}}

\newcommand{\haf}[1]{{\hat{\boldsymbol{#1}}}}
\newcommand{\ha}[1]{{\hat{#1}}}
\newcommand{\chf}[1]{{\check{\f{#1}}}}
\newcommand{\alp}{\alpha}
\newcommand{\bet}{\beta}
\newcommand{\bdg}{\beq\begin{diagram}}
\newcommand{\edg}{\end{diagram}\eeq}

\newcommand{\CNinf}{\B C_{N\infty}}
\newcommand{\CMinf}{\B C_{M\infty}}
\newcommand{\CNalp}{\B C^v_{N\alp}}
\newcommand{\CMalp}{\B C^v_{M\alp}}
\newcommand{\CNz}{\B C_{N0}}
\newcommand{\CMz}{\B C_{M0}}
\newcommand{\CNb}{\bar{\B C}_{N}}
\newcommand{\CMb}{\bar{\B C}_{M}}
\newcommand{\betgp}{\f\bet^{n-1}_{\Gam\alp} }

\newcommand{\betkp}{\f\bet^{n-1}_{K \alp } }

\newcommand{\betgps}{\f\bet^{n-1}_{\Gam\alp,s} }
\newcommand{\betkps}{\f\bet^{n-1}_{K \alp ,s} }

\usepackage{amssymb}
\usepackage{amsthm}

\newproof{prf}{Proof}[section]
\newproof{rmk}{Remark}[section]

\usepackage{lineno}
\usepackage[top=1in,bottom=1in,left=1in,right=1in]{geometry}

\journal{}%Computer Methods in Applied Mechanics and Engineering}

\begin{document}

\begin{frontmatter}

\title{An efficient displacement-based isogeometric formulation for geometrically exact viscoelastic beams}

\author[fi]{Giulio Ferri}
\author[fi]{Diego Ignesti}
\author[fi]{Enzo Marino\corref{cor1}}
\ead{enzo.marino@unifi.it}
\cortext[cor1]{Corresponding author}

\address[fi]{Department of Civil and Environmental Engineering, University of Florence\\ Via di S. Marta 3, 50139 Firenze, Italy}

\begin{abstract}
We propose a novel approach to the linear viscoelastic problem of shear-deformable geometrically exact beams. The generalized Maxwell model for one-dimensional solids is here efficiently extended to the case of arbitrarily curved beams undergoing finite displacement and rotations. 
High efficiency is achieved by combining a series of distinguishing features, that are 
i) the formulation is displacement-based, therefore no additional unknowns, other than incremental displacements and rotations, are needed for the internal variables associated with the rate-dependent material; 
ii) the governing equations are discretized in space using the isogeometric collocation method, meaning that elements integration is totally bypassed; 
iii) finite rotations are updated using the incremental rotation vector, leading to two main benefits: minimum number of rotation unknowns (the three components of the incremental rotation vector) and no singularity problems; 
iv) the same $\SO3$-consistent linearization of the governing equations and update procedures as for non-rate-dependent linear elastic material can be used;
v) a standard second-order accurate time integration scheme is made consistent with the underlying geometric structure of the kinematic problem. 
Moreover, taking full advantage of the isogeometric analysis features, the formulation permits accurately representing beams and beam structures with highly complex initial shape and topology, paving the way for a large number of potential applications in the field of architectured materials, meta-materials, morphing/programmable objects, topological optimizations, etc. 
Numerical applications are finally presented in order to demonstrate attributes and potentialities of the proposed formulation. 
\end{abstract}

\end{frontmatter}

%\linenumbers

\section{Introduction}
The interest in designing materials with predetermined mechanical properties or objects that can change shape and function over time in a programmed way is dramatically growing due to the enormous potential that can be unlocked by emerging additive manufacturing techniques. Hybrid and architectured materials \cite{Ashby&Brechet2003,Ashby2013} are examples of new ``materials'' obtainable by assembling two or more materials, possibly with a specific voids-matter distribution, such that they can exhibit attributes not owned by any (continuous distribution) of the individual materials. 
This principle is valid and applicable at micro \cite{Xu_etal2015,Shucong_etal2021,Greer_etal2023}, meso \cite{Yan_etal2016,Estrin_etal2021,Cheng_etal2023}, and macro \cite{Ge_etal2016,Ding_etal2017,Hamel_etal2019,Kirillova&Ionov2019,Boley_etal2019,Mueller_etal2022} length scales. In many cases, the inner architectures are assemblies of elements that individually can be well represented by one-dimensional deformable solids \cite{Zakharov&Pismen2019,Rafsanjani_etal2019}. Programming this type of structures represents a complex inverse problem since the design space, in addition to the properties of the single materials, comprises topology (e.g. cell shape and fibers orientation), number and connectivity of the one-dimensional elements, their (possibly locally varying) initial curvature, composition of possible multi-layered cross sections, and other parameters \cite{Florijn_etal2016,Shucong_etal2021,Wan_etal2022}. 

For such complex systems, since three-dimensional continuum-based theories would lead to prohibitive computational costs, advanced beam theories, deployed through appropriate computational formulations, can have a tremendous impact on the predictive capabilities of the numerical models. There are some essential requirements that these theories must meet: the underlying kinematic model must be able to reproduce finite displacements and rotations without any restrictions; any three-dimensional initial geometry, possibly featured by strong and local variations of curvatures, must be accurately represented; proper nonlinear material models, including rate and temperature dependent ones, need to be supported and adapted to the one-dimensional case. 

From the kinematics standpoint, starting from the fundamental works by Simo~\cite{Simo1985,SimoVu-Quoc1986} 
and Cardona \& Geradin~\cite{Cardona1988}, a number of valuable contributions have been given to the development of large deformation space beams, both to shear-deformable \cite{Dvorkin1988,Crisfield1990,Ibrahimbegovic1995,Ibrahimbegovic1995a,Ibrahimbegovic1997,Jelenic1999,Romero&Armero2002,Romero2004,Ritto-Correa&Camotim2002,Betsch&Steinmann2002,Magisano_etal2020} and Euler-Bernoulli formulations \cite{Magisano_etal2021,Vo_etal2021,Greco_etal2022,Borkovic_etal2022,Borkovic_etal2023}. 
These works are all restricted to linear elastic materials.
Inelastic material models were initially proposed in \cite{Park&Lee_1996,Gruttmann_etal2000,Battini&Pacoste_2002}. More recently, elastoplasticity has been addressed in \cite{Smriti_etal2020,Herrnboeck_etal2021,Herrnboeck_etal2023}. General three-dimensional constitutive laws are adapted to geometrically exact beams in \cite{Maqueda&Shabana_2007,Mata_etal2007,Mata_etal2008,Wackerfuss&Gruttmann_2009,Wackerfuss&Gruttmann_2011,Klinkel&Govindjee_2002}, and in \cite{Choi_etal2021,Choi_etal2022} considering also deformable directors.
Focusing on viscoelasticity, finite difference schemes are used in \cite{Lang_etal2011,Lang&Arnold2012,Linn_etal2013} with a quaternion-based parametrization of finite rotations and in \cite{Giusteri_etal2021} using the special Euclidean group for the beam kinematics. 
Finite element formulations with applications to multi-body dynamics are proposed in \cite{Zhang_etal2009,Mohamed&Shabana_2011,Bauchau&Nemani2021}.

Only recently, geometrical and material nonlinearities have been combined to specifically address the problem of programmable structures. Weeger et al.~\cite{Weeger_etal2018} propose a fully isogeometric model for geometrically exact space beams with spatially varying geometric and material properties. The possibility of modeling composite rods \cite{Ding_etal2018} (e.g., with bilayers \cite{Boley_etal2019,Zakharov_etal2020}) enables bending and twisting deformations to morph one-dimensional rods into three-dimensional shapes. In \cite{Weeger_etal2019}, although a linear elastic material is used, the advantages of a geometrically exact formulation are exploited to optimize curved beams in the sense of matching target shapes for the beam axis. In \cite{Audoly_etal2013,Lestringant_etal2020,Lestringant&Kochmann2020}, a discrete geometrically exact formulation encompassing constitutive models for incompressible viscous fluids and viscoelastic one-dimensional solids is proposed and used for simulating active and shape-morphing structures. 
Linear viscoelastic materials in a co-rotational setting are modeled in \cite{Glaesener_etal2021}, an elasto-visco-plastic material model is proposed in \cite{Weeger_etal2022}, while a two-scale approach for nonlinear hyperelastic beams is intriduced in \cite{Clezio_etal2023}. Functionally graded beams, which are also relevant for material programming, are discussed in \cite{Weeger_etal2018,Nguyen_etat2022}. 

Given the key role of time- and strain rate-dependency for active materials \cite{Bertoldi_etal2017,Dykstra_etal2019,Janbaz_etal2020,Glaesener_etal2021,Gavazzoni_etal2022} and considered the strong need for more efficient beam formulations meeting the requirements recalled above, in this paper we propose a displacement-based isogeometric collocation (IGA-C) scheme for linear viscoelastic geometrically exact beams. IGA-C \cite{Auricchio2010,Auricchio2012,Fahrendorf_etal2022} keeps the attributes of classical isogeometric analysis \cite{Hughes2005a,Cottrell2009} and, being based on the discretization of the strong form of the governing equations, completely bypasses the issues related to elements integration. The method ensures high efficiency since it requires only one evaluation point per degree of freedom, regardless of the approximation degree \cite{Schillinger2013}. IGA-C proved excellent performances for a wide range of problems \cite{Auricchio2010,Auricchio2012,Schillinger2013,Gomez2014,DeLorenzis2015,Kruse2015,Kruse2015,Gomez&DeLorenzis2016,BeiraodaVeiga2012,Auricchio2013,Kiendl2015a,KiendlAuricchioReali2017,Balduzzi_etal2017,Reali2015,Kiendl2015, KiendlMarinoDeLorenzis2017,Maurin_elal2018,Maurin_etal2018b,Evans_etal2018,Fahrendorf_etal2020,Marino_etal2020,Torre_etal2023}, including the geometrically exact beam problem \cite{Marino2016a,Weeger2017,Marino2017b,Marino2019a,Marino2019b}. 

In the present work, the generalized Maxwell model for one-dimensional solids is efficiently extended to the case of arbitrarily curved beams undergoing finite displacement and rotations. Compared with existing works, here higher efficiency is sought by combining a series of desirable and distinguishing features, in addition to the intrinsic properties of the collocation method itself. They are:
i) The formulation is displacement-based: we demonstrate that no additional unknowns with respect to the linear elastic case are needed due to the internal variables associated with the rate-dependent material. Namely, incremental displacements and rotations are the only needed unknowns; 
ii) Finite rotations are updated using the incremental rotation vector, leading to two main benefits: minimum number of rotation unknowns (i.e., the three components of the incremental rotation vector only) and no singularity problems as opposed to the case of total rotation vector; 
iii) The same $\SO3$-consistent linearization of the governing equations and update procedures as for linear elastic materials (see \cite{Marino2016a,Marino2017b,Marino2019a,Marino2019b}) can be directly used, meaning that no additional complexity in the construction the tangent matrix is introduced by the viscoelastic material; 
iv) A standard second-order accurate time integrator is made consistent with the underlying geometric structure of the kinematic problem.

In this work we also extend to the nonlinear and rate-depended case the capabilities recently developed in \cite{Ignesti_etal2023} in order to address beams and beam structures with highly complex initial curvature. 

The outline of the paper is as follows. In Section \ref{sec:bal_eqs} we review the governing equations for the quasi-static problem and the rate-dependent material model. In Section \ref{sec:time_space_disc} we introduce the three main steps to derive the formulation: time discretization, $\SO3$-consistent linearization, and space discretization. 
In Section~\ref{sec:num_app} we show through several numerical applications with increasing complexity the capabilities of the proposed formulation. Finally, in Section~\ref{sec:conclusions}, we summarize and draw the main conclusions of our work.

\section{Governing equations\label{sec:bal_eqs}}
In this section, we first introduce the strong form of the governing equations and then formulate the linear viscoelastic beam problem in a displacement-based form. 

\subsection{Strong form of the balance equations}
Let $s\mto \f c(s) \in \Rn^3$ be a space curve representing the axis of a beam, where $s\in [0, L]\sub\Rn$ is the arc length parameterization. 

The governing equations in the strong form are given in the spatial setting as follows \cite{Simo1985}
\bAl
\f{n},_{s} + \, \baf n & = 0\,,  \label{eq:spa_f} \\ 
\f{m},_{s} + \, \f{c},_{s}\car \, \f{n} +\baf m & = 0\,, \label{eq:spa_m}
\end{align}
valid for any $s\in(0,L)$ and time instant $t\in[0,T]\sub\Rn$. $L$ and $T$ are the upper bounds of the space and time domains, respectively. In the above equations, $\f{n}$ and $\f{m}$ denote the internal forces and couples per-unit length, respectively, whereas $\baf{n}$ and $\baf{m}$ are the external distributed forces and couples. 
Boundary conditions in the spatial setting write as follows
\bEq
\f \eta = \baf \eta_c \sepr{or} \f{n}  = \baf n_c \sepr{with}	s = \{0,L\}\,, t\in[0,T]\,,\label{eq:bcseta_or_n}\\
\eEq
\bEq
\f \vtht = \baf \vtht_c \sepr{or} \f m   = \baf m_c  \sepr{with} s = \{0,L\}\,, t\in[0,T]\,,\label{eq:bcstht_or_m}\\
\eEq
where $\baf\eta_c$ and $\baf\vtht_c$ are the spatial displacements and rotations imposed to any of the beam ends, $\baf n_c$ and $\baf m_c$ are the external concentrated forces and moments applied to any of the beam ends. 

Eqs.~\eqref{eq:spa_f} and \eqref{eq:spa_m} can be pulled-back in the material form as follows
\bAl
\tif K \f N +\f N,_s+\fR R\Tra \baf n =  0\,,  \label{mat_f} \\ 
\tif K \f M +\f M,_s+ \lf( \fR R\Tra \f{c},_{s}\rg) \car \, \f N  + \fR R\Tra \baf m =  0\,, \label{mat_m}
\end{align}
where $s\mto \fR R(s) \in \SO3$ is the rotation operator mapping the (rigid) rotation of the beam cross section at $s$ from the material to the current configuration. 
Note that to obtain the above material form of the balance equations, we exploited the orthogonality of $\fR R$ and the (material) beam curvature, a skew-symmetric tensor defined as $s\mto \tif K(s) \byd \fR R\Tra \fR R,_{s}\in\so3$\footnote{With the symbol $\sim$ we mark elements of $\so3$, that is the set of $3 \times 3$ skew-symmetric matrices. In this context, they are used to represent curvature matrices and infinitesimal incremental rotations. Furthermore, we recall that for any skew-symmetric matrix $\tif a\in\so3$,  $\f a = \textnormal{axial}(\tif a)$ indicates the axial vector of $\tif a$ such that $\tif a \f h= \f a \car \f h$, for any $\f h\in \Rn^3$, where $\car$ is the cross product.}. 
Moreover, $\f N = \fR R\Tra \f n$ and $\f M = \fR R\Tra \f m$ are the internal forces and couples per-unit length in the material form, respectively. Similar transformations apply to the boundary conditions. 

\subsection{The generalized Maxwell model for one-dimensional problems}
To reproduce the viscoelastic material behaviour, here we employ the generalized Maxwell model \cite{Christensen2013} directly to the one-dimensional strain measures. We first recall the material form of the strain measures \cite{Simo1985,Crisfield1999,Kapania2003}
\bEq\label{eq_strain_measures}
\f \Gam_N =  \f \Gam - \f \Gam_0 =  \fR R\Tra \f c,_s - \fR R_0\Tra \f c_{0,s}\sepr{and} \f K_M =\rm axial (\tif K - \tif K_0) = \f K - \f K_0\,, 
\eEq
where $\f \Gam = \fR R\Tra\f c,_s - [1,\, 0,\, 0]\Tra$ and $\f \Gam_0 = \fR R_0\Tra\f c_{0,s} - [1,\, 0,\, 0]\Tra$. 
$\f \Gam_N\in\Rn^3$ describes the axial and shear strains, whereas $\f K_M\in\Rn^3$ describes the bending and torsional strains; 
$\tif K_0 = \fR R_0\Tra \fR R_{0,s}\in\so3$ is the beam initial curvature (skew-symmetric tensor) in the material form; $ \f c_0$ represents the beam axis in the initial configuration;
$\fR R_0\in\SO3$ is the rotation operator that expresses the rotation of the beam cross section in the initial configuration. 

Assuming a rheological model made of \textit{m} spring-dashpot elements, the material form of the total internal forces and couples can be written as
\bEq\label{stress_viscoel}
\f N = \f N_{\infty} +\sum_{\alp=1}^{m}\f N_{\alp}\sepr{and} \f M = \f M_{\infty} +\sum_{\alp=1}^{m}\f M_{\alp}\,,
\eEq
where $\f N_{\infty}$ and $\f M_{\infty}$ are the long-term elastic internal forces and couples, whereas $\f N_{\alp}$ and $\f M_{\alp}$ are the viscous contributions related to the $\alp$th Maxwell element. Assuming the existence of a one-dimensional dissipative potential, additively made of a contribution related to axial and shear strains and a contribution related to torsional and bending strains (see also \cite{Lang&Arnold2012,Weeger_etal2022}), the dissipative internal stresses are given by 
\bEq
\f N_{\alp}=\B H^v_{N\alp}\dot{\f\Gam}_{N\alp} \sepr{and} \f M_{\alp}=\B H^v_{M\alp}\dot{\f K }_{M\alp}\,,\label{if_visco}
\eEq
where $\dot{\f\Gam}_{N\alp}$ and $\dot{\f K }_{M\alp}$ are the viscous strain rate vectors and 
$\B H^v_{N\alp}$ and $\B H^v_{M\alp}$ are diagonal viscosity matrices associated with the $\alp$th Maxwell element. 
Introducing the relaxation times, $\tau_\alp$ with $\alp = 1,\ldots,m$, that without loss of generality are assumed to be the same for all the strain measures, the rates of viscous strains for the $\alp$th Maxwell element can be expressed through the evolutionary equations as follows
\bEq\label{evol_f_m}
\dot{\f\Gam}_{N\alp}=\fr {1} {\tau_\alp}(\f\Gam_N-\f\Gam_{N\alp})
\sepr{and}
\dot{\f K}_{M\alp}=\fr {1} {\tau_\alp}(\f K_M-\f {K}_{M\alp})\,. 
\eEq

Finally, assuming that the elastic response is linear, the total material internal forces and couples become 
\bAl
\f N &= \CNinf\f\Gam_N +\sum_{\alp=1}^{m}\CNalp(\f\Gam_N-\f\Gam_{N\alp})\,,\label{N_constitutive}\\
\f M &= \CMinf\f{K}_M +\sum_{\alp=1}^{m}\CMalp(\f{K}_M-\f{K}_{M\alp})\,,\label{M_constitutive}
\end{align}
with 
\beq\label{Cinf}
\CNinf=\mathrm{diag}(E_{\infty}A,G_{\infty}A_2,G_{\infty}A_3) \sepr{and} \CMinf=\mathrm{diag}(G_{\infty}J_t,E_{\infty}J_2,E_{\infty}J_3)\,,
\eeq
\beq\label{CNalp}
\CNalp=\B H^v_{N\alp}/\tau_\alp = \mathrm{diag}(E_{\alp}A,G_{\alp}A_2,G_{\alp}A_3)
\sepr{and}
\CMalp= \B H^v_{M\alp}/\tau_\alp =\mathrm{diag}(G_{\alp}J_t,E_{\alp}J_2,E_{\alp}J_3)\,,
\eeq
where $E_{\infty}$ and $G_{\infty} = E_{\infty}/2(1+\nu)$ are the long term Young and shear moduli, 
whereas $E_{\alp}$ and $G_{\alp} = E_{\alp}/2(1+\nu)$ are the Young and shear moduli associated with the $\alp$th Maxwell element. 
We assume a constant Poisson ratio for the material.

\subsection{Displacement-based strong form of the governing equations}
By substituting Eqs.~\eqref{N_constitutive} and \eqref{M_constitutive} into \eqref{mat_f} and \eqref{mat_m},
we obtain the governing equations expressed in terms of kinematic quantities only, i.e. total and viscous strain measures, as follows
\begin{gather}
\tif K \CNz\f\Gam_N - \tif K\sum_{\alp=1}^{m} \CNalp \f\Gam_{N\alp} + \CNz\f\Gam_{N,s} - \sum_{\alp=1}^{m}\CNalp\f\Gam_{N\alp,s}+\fR R\Tra \baf n = \f 0\,,\label{eq:mat_f_strain}\\
\tif K \CMz\f{K}_M-\tif K\sum_{\alp=1}^{m}\CMalp\f{K}_{M\alp}+\CMz\f K_{M,s} - \sum_{\alp=1}^{m}\CMalp\f K_{M\alp,s} + \fR R\Tra\f{c},_{s}\car{\CNz\f\Gam_N}+ \nonumber \\
 - \fR R\Tra\f{c},_{s}\car{\sum_{\alp=1}^{m}\CNalp\f\Gam_{N\alp}}+\fR R\Tra \baf m =  \f 0\,.\label{eq:mat_m_strain}
\end{gather}

Similarly, the Neumann boundary conditions become
\bEq\label{eq:mat_bcsf_strain} 
\CNz\f\Gam_N-\sum_{\alp=1}^{m}\CNalp\f\Gam_{N\alp}=\fR R\Tra \baf n_c\,,
\eEq
\bEq\label{eq:mat_bcsmm_strain} 
\CMz\f{K}_M-\sum_{\alp=1}^{m}\CMalp\f{K}_{M\alp}=\fR R\Tra \baf m_c\,,
\eEq
where $\CNz = \CNinf + \sum_{\alp=1}^m \CNalp$ and $\CMz =  \CMinf + \sum_{\alp=1}^m \CMalp$ are the instantaneous elasticity tensors.

\section{Discretization and consistent linearization of the governing equations\label{sec:time_space_disc}}
In this section, after introducing the time discretization, we present the consistent linearization and space discretization of the governing equations. 

\subsection{Time discretization and integration scheme}
By using the standard trapezoidal rule for the time integration, at time $t^n = (n-1)dt$, where $dt$ is the time step span and $n=1,\ldots$ is the time step counter, 
the time discretized viscous strain measures $\f\Gam^n_{N\alp}$ and $\f{K}^n_{N\alp}$ associated with the $\alp$th Maxwell element can be expressed as
\bEq\label{g_trapz}
\f\Gam^{n}_{N\alp}=\fr {dt} {(2\tau_\alp+dt)}\f\Gam^{n}_N+\betgp\,,
\eEq
\bEq\label{k_trapz}
\f {K}^{n}_{N\alp}=\fr {dt} {(2\tau_\alp+dt)}\f{K}^{n}_M+\betkp\,,
\eEq
where the terms $\betgp$ and $\betkp$ are defined as follows 
\bEq\label{betg}
\betgp=\fr {dt} {(2\tau_\alp+dt)}\f\Gam^{n-1}_N+\fr{2\tau_\alp-dt} {(2\tau_\alp+dt)}\f\Gam^{n-1}_{N\alp}\,,
\eEq
\bEq\label{betk}
\betkp=\fr {dt} {(2\tau_\alp+dt)}\f{K}^{n-1}_M+\fr{2\tau_\alp-dt} {(2\tau_\alp+dt)}\f{K}^{n-1}_{M\alp}\,.
\eEq

Importantly, we remark that the above terms are computed using only quantities known form the previous time step, therefore they do not need to be iteratively updated during the Newton-Raphson algorithm.

Substituting Eqs.~\eqref{g_trapz} and~\eqref{k_trapz} into the time discretized form of Eqs.~\eqref{eq:mat_f_strain} and \eqref{eq:mat_m_strain}, the nonlinear governing equations can be expressed in terms of total strains (and other quantities known from the previous time step) as follows
\begin{gather} 
\tif K^n[\CNz-\sum_{\alp=1}^{m}\CNalp\fr {dt} {(2\tau_\alp+dt)}]\f\Gam^{n}_N - \tif K^{n}\sum_{\alp=1}^{m}\CNalp\betgp + \nonumber\\
+[\CNz-\sum_{\alp=1}^{m}\CNalp\fr {dt} {(2\tau_\alp+dt)}]\f\Gam^{n}_{N,s}-\sum_{\alp=1}^{m}\CNalp\betgps+\fR {R\Tra}^n \baf n^{n} =  \f 0\,,\label{eq:mat_f_def}
\end{gather}

\begin{gather} 
\tif K^{n}[\CMz-\sum_{\alp=1}^{m}\CMalp\fr {dt} {(2\tau_\alp+dt)}]\f{K}^{n}_M+[\CMz-\sum_{\alp=1}^{m}\CMalp\fr {dt} {(2\tau_\alp+dt)}]\f{K}^{n}_{M,s}+ \nonumber\\
-\tif K^{n}\sum_{\alp=1}^{m}\CMalp\betkp-\sum_{\alp=1}^{m}\CMalp\betkps+\fR R\Tra{^n}\f{c}^{n},_{s}\car{[\CNz-\sum_{\alp=1}^{m}\CNalp\fr {dt} {(2\tau_\alp+dt)}]\f\Gam^{n}_N}+ \nonumber\\
-\fR R\Tra{^n}\f{c}^{n},_{s}\car{\sum_{\alp=1}^{m}\CNalp\betgp}+\fR R\Tra{^n} \baf m^{n} =  \f 0\,.\label{eq:mat_m_def}
\end{gather}

The same can be done for the Neumann boundary conditions 
\bEq\label{eq:mat_bcsf_def} 
 [\CNz-\sum_{\alp=1}^{m}\CNalp\fr {dt} {(2\tau_\alp+dt)}]\f\Gam^{n}_N=\sum_{\alp=1}^{m}\CNalp\betgp+\fR {R\Tra}^{n} \baf n^{n}_c\,,
\eEq
\bEq\label{eq:mat_bcsmm_def} 
[\CMz-\sum_{\alp=1}^{m}\CMalp\fr {dt} {(2\tau_\alp+dt)}]\f{K}^{n}_M=\sum_{\alp=1}^{m}\CMalp\betkp+\fR R\Tra{^n} \baf m^{n}_c\,.
\eEq

\subsection{$\SO3$-consistent linearization of the displacement-based governing equations}
A remarkable advantage of the present formulation is that the same $\SO3$-consistent linearization rules derived for static or dynamic formulations with non-rate-dependent material, see e.g. \cite{Marino2019b}, can be directly used. 
This is possible since, as noticed above, the viscous components can be expressed by terms depending on the total strains (likewise the non-rate-dependent case) and terms ($\betgp$ and $\betkp$) which are known from the previous time step and therefore do not substantially affect the construction of the tangent operator.

To preserve $\SO3$-consistency, we only recall some fundamental geometrical aspects. Given a curve in the configuration manifold $\veps \mto\f\gam(\veps) = (\f c_\veps,\fR R_\veps)$, with $\veps \in \Rn$, defined by $\f c_\veps = \f c + \veps \del \f \eta$ (standard translation on $\Rn^3$) and $\fR R_\veps = \fR R \exp(\veps \del \tif \vTht)$ (left translation on $\SO3)$, the linearization is based on the construction of the (material) tangent space at $(\f c,\fR R)$ obtained by $(d \f \gam/d \veps)_{\veps = 0}$ such that $\f\gam(0) = (\f c,\fR R)$. From the kinematic point of view, $\del \f\eta\in\Rn^3$ represents an incremental displacement superimposed to the current configuration of the centroid line $\f c$; whereas $\del \tif \vTht \in \so3$ represents an incremental rotation superimposed to the rotation $\fR R$. Note that in the construction of the curve $\f \gam$ we used the exponential map $\exp :\so3 \to \SO3$ which maps the line $\veps\tif\vTht$  at $\so3$ onto the one parameter subgroup $\exp(\veps\tif\vTht) \in\SO3$ \cite{Choquet-Bruhat1996}. Note that for $\SO3$ the exponential map is expressed by an exact (Rodrigues) formula.

By using the linearizion rules given in \cite{Marino2019b}, the linearized version of Eqs.~\eqref{eq:mat_f_def}-\eqref{eq:mat_bcsmm_def} is 
\begin{gather}
[\CNb\wti{(\hat{\fR R}{\Tra{^n}}\hat{\f{c}}^{n}_{,s})}-\wti{(\CNb\hat{\f\Gam}^{n}_N)}+\sum_{\alp=1}^{m}\wti{(\CNalp\betgp)}]\del\f\vTht,_s^{n}+\nonumber\\
+[\ha{\tif K^{n}}\CNb\wti{(\hat{\fR R}{\Tra{^n}}\hat{\f{c}}^{n}_{,s})}-\wti{(\CNb\hat{\f\Gam}^{n}_N)}\ha{\tif K^{n}}+\ha{\tif K^{n}}\sum_{\alp=1}^{m}\wti{(\CNalp\betgp)}-\wti{(\ha{\tif K^{n}}\sum_{\alp=1}^{m}{\CNalp\betgp})}+\nonumber\\
+\CNb\wti{(\hat{\fR R}{\Tra{^n}}\hat{\f{c}}^{n}_{,ss})}-\CNb\wti{(\tif K^{n}\hat{\fR R}{\Tra{^n}}\hat{\f{c}}^{n}_{,s})}+\wti{(\hat{\fR R}{\Tra{^n}}{\baf n}^{n})}]\del\f\vTht^{n}+\CNb\hat{\fR R}{\Tra{^n}}\del\f\eta_{,ss}^n+\nonumber\\
+[\tif K^{n}\CNb\hat{\fR R}\Tra{^n}-\CNb\tif K^{n}\hat{\fR R}{\Tra{^n}}]\del\boldsymbol\eta_{,s}^n+\hfR F^n = \f 0\,,\label{traslin}
\end{gather}

\begin{gather}
\CMb\del\f\vTht,_{ss}^{n}+[\CMb\ha{\tif K^{n}}+\ha{\tif K^{n}}\CMb-\wti{(\CMb\hat{\f{K}}^{n}_M)}+\sum_{\alp=1}^{m}\wti{(\CMalp\betkp)}]\del\f\vTht,_s^{n}+\nonumber\\
+[\ha{\tif K^{n}}\CMb\ha{\tif K^{n}}-\wti{(\CMb\hat{\f{K}}^{n}_M)}\ha{\tif K^{n}}+\ha{\tif K^{n}}\sum_{\alp=1}^{m}\wti{(\CMalp\betkp)}-\wti{(\ha{\tif K^{n}}\sum_{\alp=1}^{m}{\CMalp\betkp})}+\nonumber\\
+\CMb\ha{\tif K_{,s}^{n}}+[\wti{(\hat{\fR R}\Tra{^n}\hat{\f{c}}^{n}_{,s})}\CNb-\wti{(\CNb\hat{\f\Gam}^{n}_N)}+\sum_{\alp=1}^{m}\wti{(\CNalp\betgp)}]\wti{(\hat{\fR R}\Tra{^n}\hat{\f{c}}^{n}_{,s})}+\wti{(\hat{\fR R}\Tra{^n}\baf m^{n})}]\del\f\vTht^{n}+\nonumber\\
+[\wti{(\hat{\fR R}\Tra{^n}\haf c,_s^n})\CNb-\wti{(\CNb\hat{\f\Gam}^{n}_N)}+\sum_{\alp=1}^{m}\wti{(\CNalp\betgp)}]\hat{\fR R}\Tra{^n}\del\boldsymbol\eta^{n}_{,s}+\hfR T^n =\f 0\,,\label{rotlin}
\end{gather}
where we have defined
\begin{gather} 
\hfR F^n =\ha{\tif K^{n}}\CNb\hat{\f\Gam}^{n}_N-\ha{\tif K^{n}}\sum_{\alp=1}^{m}\CNalp\betgp+\CNb\hat{\f\Gam}^{n}_{N,s}-\sum_{\alp=1}^{m}\CNalp\betgps+\hat{\fR R}\Tra{^n} \baf n^{n}\,, \label{F_trasl} \\
\hfR T^n = \ha{\tif K^{n}}\CMb\hat{\f{K}}^{n}_M - \ha{\tif K^{n}}\sum_{\alp=1}^{m}\CMalp\betkp+\CMb\hat{\f{K}}^{n}_{M,s}-\sum_{\alp=1}^{m}\CMalp\betkps+\nonumber\\
+\hat{\fR R}\Tra{^n}\haf c^n,_s \car{[\CNb\hat{\f\Gam}^{n}_N-\sum_{\alp=1}^{m}\CNalp\betgp]}+\hat{\fR R}\Tra{^n} \baf m^{n}\,,
\end{gather}
and we have set $\CNb=[\CNz-\sum_{\alp=1}^{m}\CNalp\fr {dt} {(2\tau_\alp+dt)}]$ and $\CMb=[\CMz-\sum_{\alp=1}^{m}\CMalp\fr {dt} {(2\tau_\alp+dt)}]$.

Similarly, for the boundary conditions we have 
\bEq 
 [\CNb\wti{(\hat{\fR R}\Tra{^n}\hat{\f{c}}^{n},_s)}-\wti{(\hat{\fR R}\Tra{^n}{\baf n}^{n}_c)}]\del\f\vTht^{n}+\CNb\hat{\fR R}\Tra{^n}\del\f\eta^{n}_{,s}=\sum_{\alp=1}^{m}\CNalp\betgp-\CNb\hat{\f\Gam}^{n}_N+\hat{\fR R}\Tra{^n} \baf n^{n}_c\,,\label{bcsf_lin}
\eEq
and 
\bEq
\CMb\del\f\vTht,_s^{n}+[\CMb\ha{\tif K^{n}}-\wti{(\hat{\fR R}\Tra{^n}{\baf m}^{n}_c)}]\del\f\vTht^{n}=\sum_{\alp=1}^{m}\CMalp\betkp-\CMb\hat{\f{K}}^{n}_M+\hat{\fR R}\Tra{^n} \baf m^{n}_c\,. \label{bcsmm_lin} 
\eEq

In the above linearized equations, with the symbol $\hat{(\cdot)}$ we denote any quantity evaluated at the time $t^n$ around which the linearization takes place.
Note that the linearized equations has just more terms, but do not add any significant complexity to the standard liner elastic rate-independent IGA-C formulations \cite{Marino2016a,Marino2017b,Marino2019a,Marino2019b}. 

\subsection{Space discretization}
The linearized governing equations, written in terms of the unknown fields $\del\f\vTht^{n}$ and $\del\f\eta^n$ at time  $t^{n}$, are spatially discretized by using NURBS basis functions $R_{j,p}$ with $j = 1,\ldots,\rm n$ of degree $p$. In the isoparametric formulation, we use the same basis functions to represent the beam centroid curve. Thus, we have
\bAl
\del\f\vTht^{n} (u) &= \sum_{j= 1}^{\rm n} R_{j,p}(u) \del \chf\vTht^{n}_j\sepr{with} u\in [0,\,1]\,,\label{eq:Thtu}\\
\del\f\eta^{n} (u) &= \sum_{j= 1}^{\rm n} R_{j,p}(u)\del \chf\eta^{n}_j\sepr{with} u\in [0,\,1]\,,\label{eq:etau}\\
\f c^{n} (u) & = \sum_{j= 1}^{\rm n} R_{j,p}(u) \chf p^{n}_j\sepr{with} u\in [0,\,1]\,,\label{eq:cu}
\end{align}
where $\del \chf\vTht^{n}_j$ and $\del\chf\eta^{n}_j$ are the primal ($2\car3\car \rm n$) unknowns, namely the $j$th incremental control rotation and translation, respectively; $\chf p^{n}_j$ is the $j$th control point defining the beam centroid curve. 
Eqs. \eqref{eq:Thtu} and \eqref{eq:etau} are substituted into \eqref{traslin} and \eqref{rotlin} and in the boundary conditions \eqref{bcsf_lin} and \eqref{bcsmm_lin}, where the differentiations must be properly done considering the Jacobian relating the parametric $u$ and arc length $s$ coordinate systems.
A square system is finally built by collocating the field equations at the internal $\rm n-2$ collocation points and the Dirichlet or Neumann boundary conditions at the boundaries $u = 0$ and $u = 1$. Standard Greville abscissae \cite{Auricchio2010} are chosen as collocation points. 
At a given Newton-Raphson iteration, once solved the linear system for the primal variables, we follow the same updating procedure discussed in \cite{Marino2019b} and make use Eqs. \eqref{g_trapz}--\eqref {betk} to update the viscous terms which are only evaluated at the collocation points.

\section{Numerical results\label{sec:num_app}}
In this section we report a series of numerical applications aimed at testing all the attributes of the proposed model. 
We start with the roll-up of a straight cantilever beam under a concentrated couple at the free-end. For the same beam, we also present the results under a concentrated tip force with a complex variation in time.
Then, a circular arch with an out-of plane tip load is considered. In this case, in order to check the rates of convergence in space, we keep the problem linear to use the exact solution as a reference. 
For all the above mentioned tests, the material properties are intentionally set to emphasize the viscoelastic behavior and to check the robustness of the proposed method.
Subsequently, the capability to model curved beam problems in the nonlinear regime is tested starting from a circular arch and then addressing the Spivak \cite{Ignesti_etal2023} and the Lissajous \cite{Marino_etal2020} beams, both featuring very challenging curvatures.
Finally, the full potentialities of the proposed formulation are shown by simulating the complex deformations of planar and cylindrical nets with cells made of curved beam elements.

\subsection{Roll-up of a cantilever beam}
The beam is placed along the $x_1$-axis of a fixed Cartesian reference system and is subjected to a concentrated couple at the free end with axis $x_2$. The length of the beam is \SI{1}{m}. The cross section is a square of side \SI{0.1}{m}. One Maxwell element is employed, i.e., $m = 1$. The Young modulus $E_1$ and relaxation time $\tau_1$ are reported in Table \ref{tab:rollup}. The shear modulus is set as $G_1 = {E_1}/{2(1+\nu)}$, where, for this specific test, $\nu=0$.
The relaxation time is $\tau_{1}=\SI{5000}{s}$ to study the convergence of the instantaneous response ($t\approx \SI{0}{s}$), whereas $\tau_{1}=\SI{0.05}{s}$ to study the convergence of the long term ($t\approx{t_{\infty}}$) response. In the former case the time step is $\SI{1E-4}{s}$, in the latter $\SI{5}{s}$.

\begin{table}\label{tab:rollup}
\footnotesize
\centering
\begin{tabular}{lcccc}
$E$ $[\SI{}{N/m^2}]$ 	& $\tau(t\approx{\SI{0}{s}})$ $[\SI{}{s}]$ & $\tau(t\approx{t_{\infty}})$  $[\SI{}{s}]$\\
\hline
40 & - & - \\
 530  	&5000 & 0.05\\
\end{tabular}
\caption{Material mechanical properties for the roll-up case.}\label{tab:rollup}
\end{table}

\begin{figure}
\centering
	\subfigure[Convergence plots for $p=4, 6, 8$ of $err_{tip}$ vs. number of collocation points.\label{fig:roll0l2}]{\includegraphics[width=0.48\textwidth]{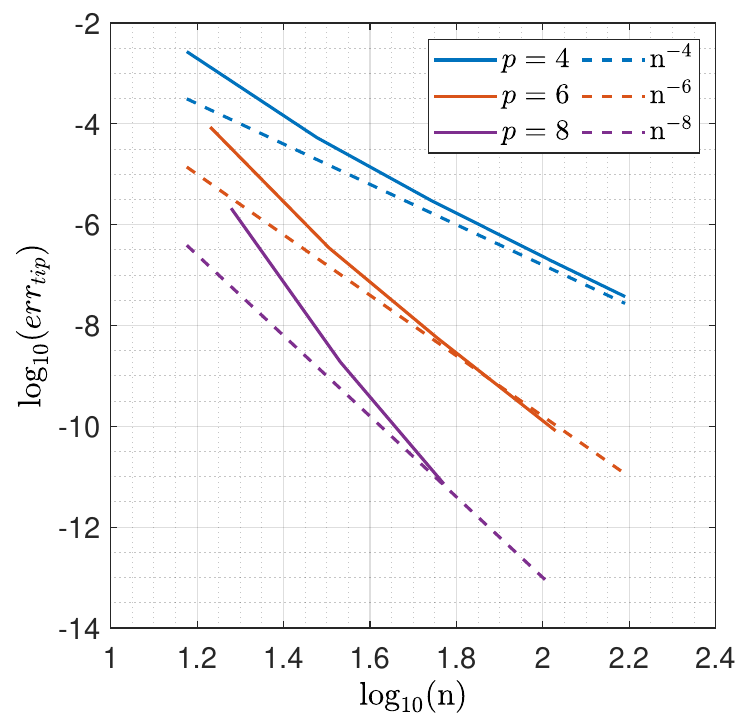}}\hspace{0.5cm}
  	\subfigure[Convergence plots for $p=4, 6, 8$ of $err_{radius}$ vs. number of collocation points.\label{fig:roll0rad}]{\includegraphics[width=0.48\textwidth]{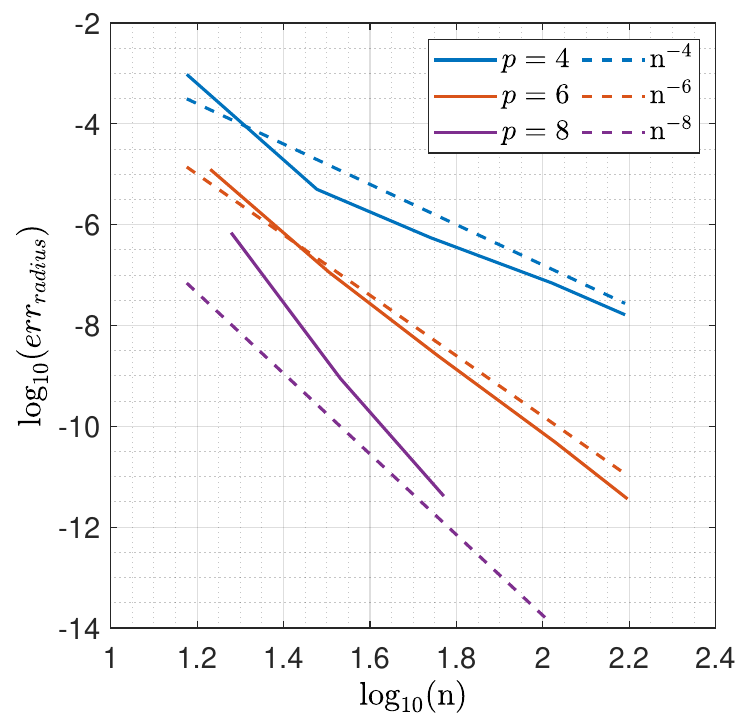}}
\caption{\label{fig:rollup0} Instantaneous roll-up of a viscoelastic cantilever beam: errors convergence.}
\end{figure}

\begin{figure}
\centering
	\subfigure[Convergence plots for $p=4, 6, 8$ of $err_{tip}$ vs. number of collocation points.\label{fig:rollinfl2}]{\includegraphics[width=0.48\textwidth]{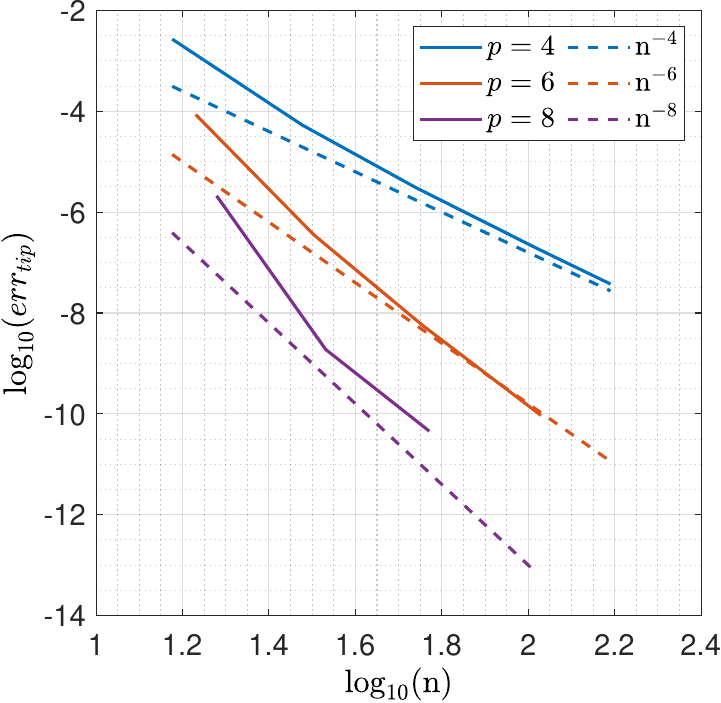}}\hspace{0.5cm}
 	\subfigure[Convergence plots for $p=4, 6, 8$ of $err_{radius}$ vs. number of collocation points.\label{fig:rollinfrad}]{\includegraphics[width=0.48\textwidth]{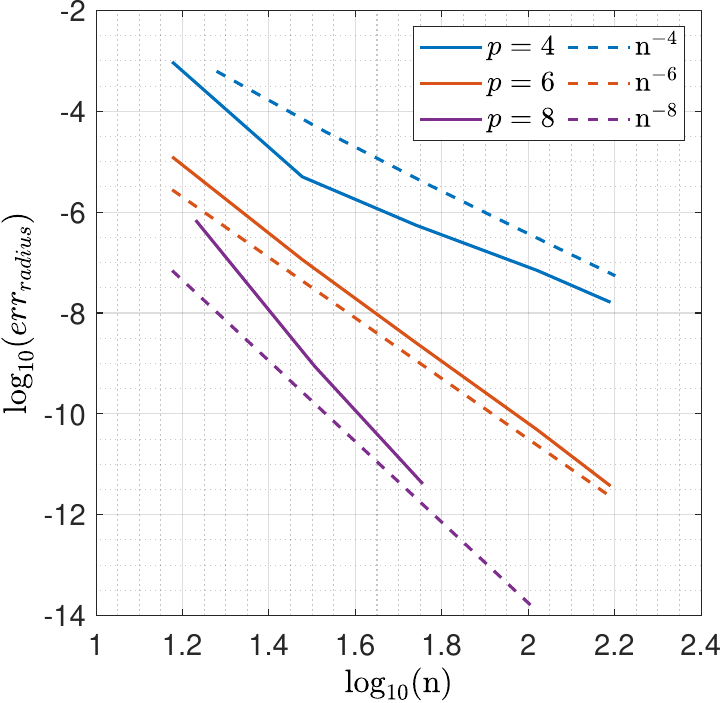}}
\caption{\label{fig:rollupinf}  Long term roll-up of a viscoelastic cantilever beam: errors convergence.}
\end{figure}

We apply a couple of magnitudes $M_{0}={2\pi(E_{\infty}+E_{1})J_{2}}/{L}$ and $M_{\infty}={2\pi E_{\infty}J_{2}}/{L}$, such that the beam deforms into a full closed circle of radius $L/2\pi$ at both times $t\approx\SI{0}{s}$ and $t\approx t_{\infty}$.
We assess the quality of the convergence rates by comparing the computed position of the beam tip, $err_{tip} = || \f{p}_{\rm n}^{h} ||_{L_2}$, with the exact value, namely the clamped end located at the origin of the reference system.
Moreover, to have also a global measure of the error, in addition to $err_{tip}$, we compare also the radius of the closed circle computed numerically, $r^{h}$, with the analytical one. $r^{h}$ is obtained selecting three points uniformly spaced on the deformed beam axis. The relative error is then computed as $err_{radius} = 2\pi |r^{h} - L/2\pi |/L$. These two errors are plotted vs. the number of collocation points and for different degrees for both instantaneous (Figure~\ref{fig:rollup0}) and long term (Figure~\ref{fig:rollupinf}) responses. 

It is observed that the two error measures exhibit very good spatial convergence rates for all the degrees considered both at instantaneous and infinity times. It is noted that the accuracy cannot be smaller than $\SI{1E-10}{}$ due to the tolerance set in the Newton-Raphson algorithm. Odd degrees are not considered in this study since, as it is well known, in collocation they do not normally improve the rates compared to the smaller even degree.

\subsection{Cantilever beam subjected to a tip force}
\begin{figure}
\centering
\subfigure[Load intensity in $x_3$ direction.\label{fig:mensolaF}]{\includegraphics[width=0.27\textwidth]{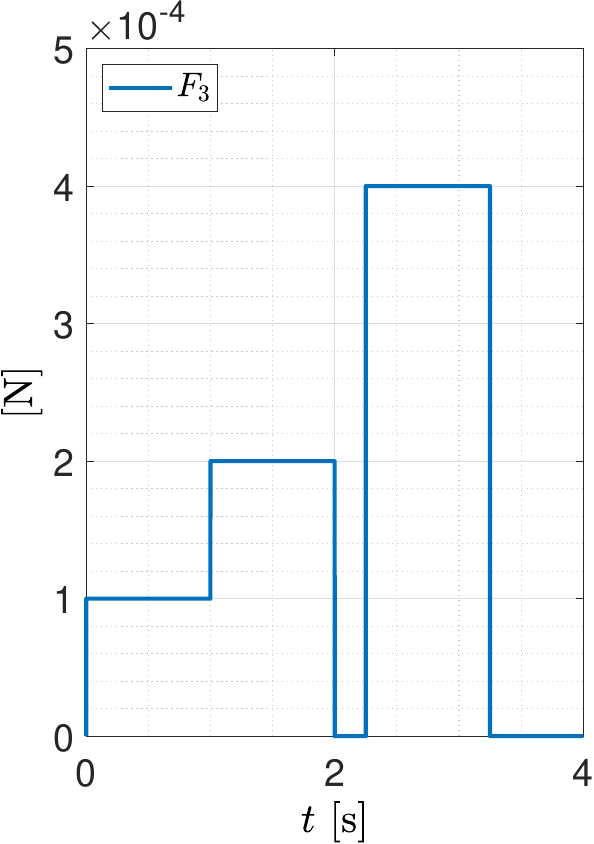}}\hspace{0.25cm}
\subfigure[Tip displacement along the $x_1$-axis.\label{fig:mensolau1}]{\includegraphics[width=0.3\textwidth]{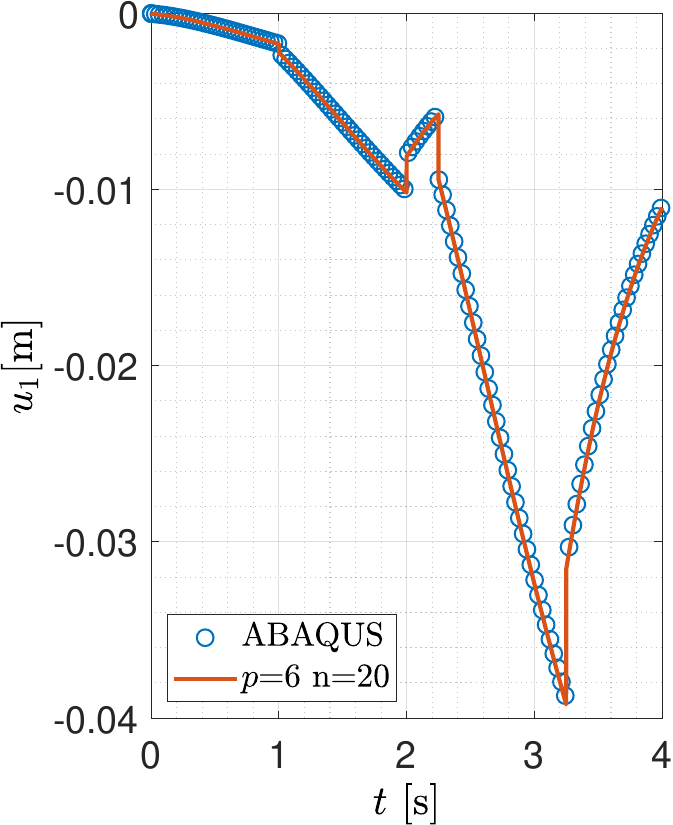}}\hspace{0.25cm}
\subfigure[Tip displacement along the $x_3$-axis.\label{fig:mensolau3}]{\includegraphics[width=0.3\textwidth]{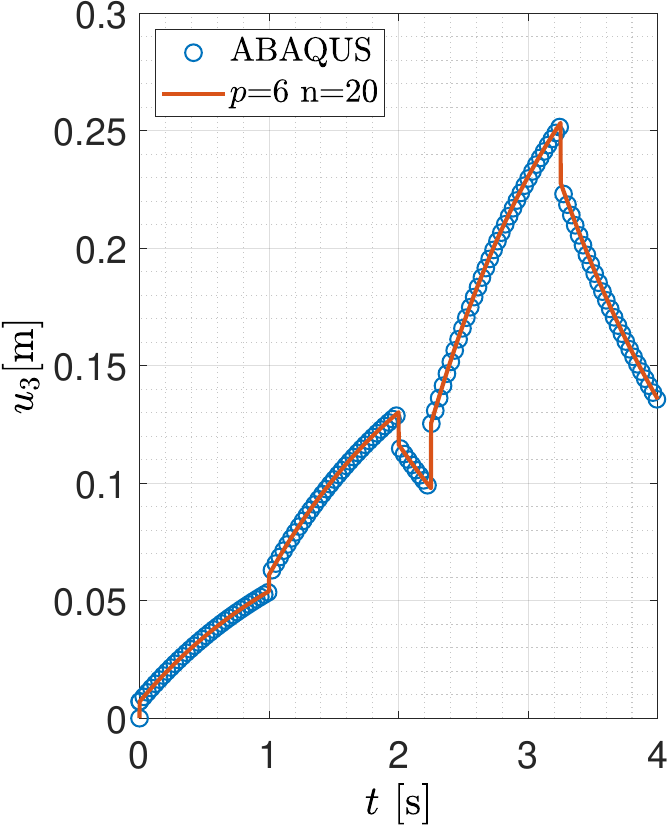}}
\caption{\label{fig:tipload_c} Cantilever beam subjected to tip load: loads and tip displacements over time.}
\end{figure}

The same cantilever beam analyzed in the previous section is here loaded with a tip force in the $x_3$ direction with the complex pattern shown in Figure~\ref{fig:mensolaF}. 
The relaxation time is $\SI{0.1}{s}$. A time step of $\SI{1E-3}{s}$ is used for a total simulation time of $\SI{4}{s}$. The time history of the displacements in the $x_1$ and $x_2$ directions are shown in Figures~\ref{fig:mensolau1} and~\ref{fig:mensolau3}, respectively, where a comparison with Abaqus \cite{Abaqus} is included. In Abaqus, 100 B31 shear-deformable beam elements have been used.
As it can be noticed from Figure~\ref{fig:mensolau3}, the beam undergoes large deflections. The viscous deformations are rather significant since, under repeated force steps, the deflection keeps increasing. For example, under a constant force of $\SI{4e-4}{N}$, from $\SI{2.25}{s}$ to $\SI{3.25}{s}$, the vertical tip displacement increases from $\SI{0.1255}{m}$ to $\SI{0.2525}{m}$. An excellent agreement with Abaqus is observed. 

\subsection{Circular arch subjected to an out-of plane tip force}
The \ang{90} circular arch here considered has a radius $R = \SI{1}{m}$ and a square cross section of side $\SI{0.1}{m}$. 
The beam lies in the $(x_1,x_2)$ plane, it is clamped at one end and subjected to a constant force along the $x_3$-direction applied at its free end. 
The beam centroid is exactly reconstructed using 3 control points and NURBS basis functions of degree 2. $k$-refinement is then applied. 
The same material properties for the long term response of the above roll-up case are here employed; the time step is $\SI{0.1}{s}$. 
The instantaneous and long term analytical solutions \cite{Auricchio2013} for this benchmark test, under the hypothesis of small displacements, can be calculated as follows 
\begin{align}
v_0      & =\frac{\pi F_{3}R}{2G_{0}A_3}+\frac{F_{3}R^{3}}{G_{0}J_t}(\frac{3}{4}\pi-2)+\frac{\pi F_{3}R^3}{4E_{0}J_2}\,, \label{tiploadcircle_0}\\
v_\infty & =\frac{\pi F_{3}R}{2G_{\infty}A_3}+\frac{F_{3}R^{3}}{G_{\infty}J_t}(\frac{3}{4}\pi-2)+\frac{\pi F_{3}R^3}{4E_{\infty}J_2}\,,\label{tiploadcircle_inf}
\end{align}
where $G_{0}=G_{\infty}+G_{1}$, $E_{0}=E_{\infty}+E_{1}$, 
$A_3$ is the beam cross section multiplied by the shear correction factor (equal to $\SI{5/6}{}$), whereas $J_t$ and $J_2$ are the torsional and the bending moment of inertia, respectively. The out-of-plane force is $F_{3} = \SI{1e-8}{N}$. 

\begin{figure}
\centering
	\subfigure[Convergence plots for $p=4,6,8$ of $err_{tip}$ for $t\approx \SI{0}{s}$ vs. number of collocation points.\label{fig:conv_c0}]{\includegraphics[width=0.48\textwidth]{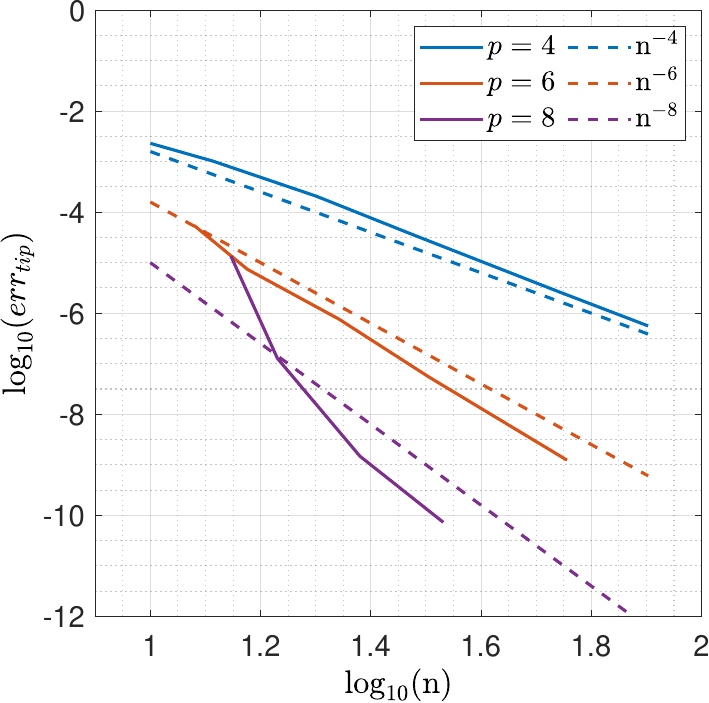}}\hspace{0.5cm}
 	\subfigure[Convergence plots for $p=4, 6, 8$ of $err_{tip}$ for $t\approx{t_\infty}$ vs. number of collocation points.\label{fig:conv_cinf}]{\includegraphics[,width=0.48\textwidth]{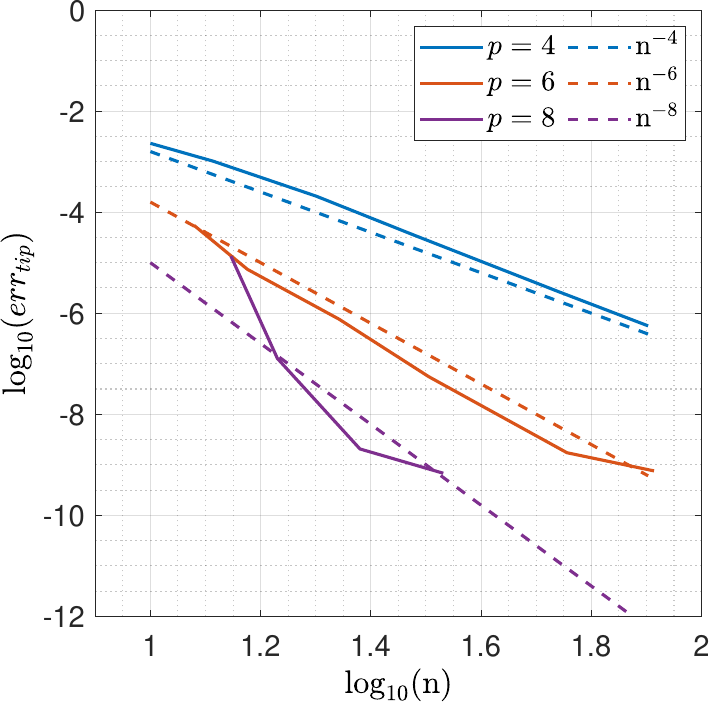}}
\caption{\label{fig:conv_c} \ang{90} circular arch subjected to an out-of plane tip force: errors convergence.}
\end{figure}

The relative errors at the beam tip for the instantaneous and the long term responses are shown in Figures~\ref{fig:conv_c0} and \ref{fig:conv_cinf}, respectively. 
Very good convergence rates are observed, although in the long term case apparently the error reaches about \SI{1E-9}{}, whereas for the instantaneous case it reaches \SI{1E-10}{} that is the tolerance limit in the Newton-Raphson algorithm.

\subsection{Circular arch under complex tip load patterns}
\begin{figure}
\centering
\begin{overpic}[width=0.7\textwidth]{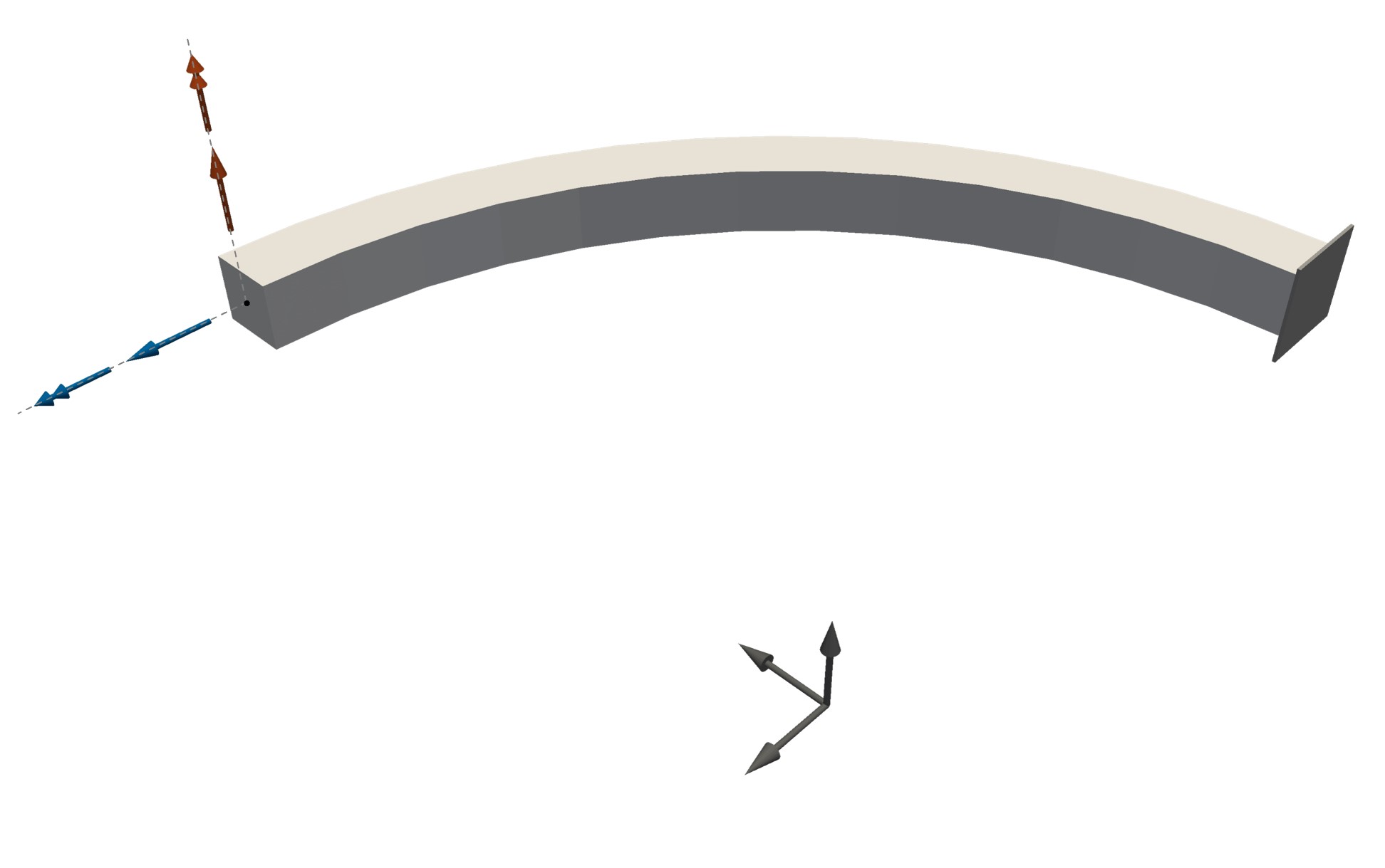}
\put(490,150){$x_1$}
\put(490,50){$x_2$}
\put(612,175){$x_3$}
\put(50,295){$M_2$}
\put(120,325){$F_2$}
\put(75,530){$M_3$}
\put(100,450){$F_3$}
\end{overpic}
\caption{\ang{90} circular arch clamped at on end and subjected to tip forces and couples along $x_2$-axis (blue) and $x_3$-axis (orange).\label{fig:cerchio_pattern}}
\end{figure}

\begin{figure}
\centering
	\subfigure[Tip forces.\label{fig:lpF}]{\includegraphics[width=0.48\textwidth]{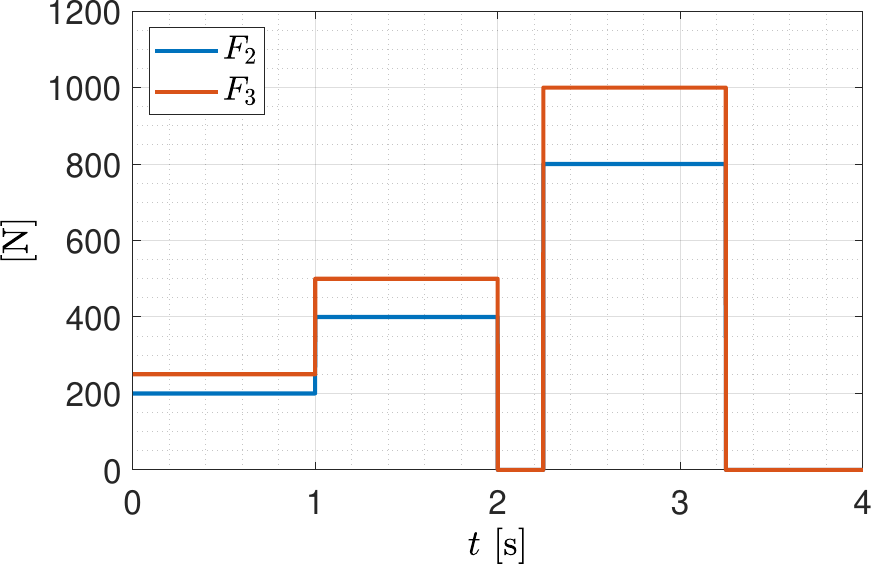}}\hspace{0.2cm}
 	\subfigure[Tip couples.\label{fig:lpM}]{\includegraphics[width=0.48\textwidth]{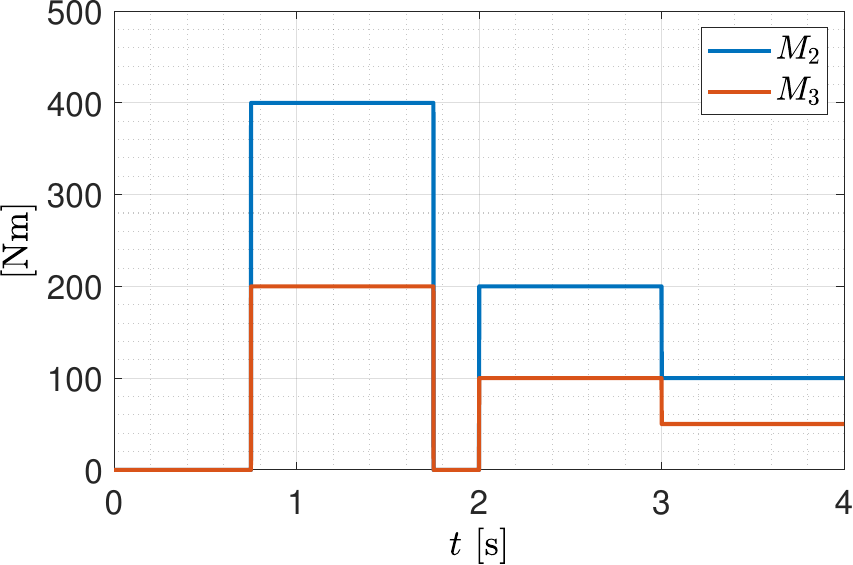}}
\caption{\label{fig:cerchio_FM} \ang{90} circular arch: tip loads time histories.}
\end{figure}

The same circular arch of the previous section is here subjected to a complex tip load pattern. As shown in Figure~\ref{fig:cerchio_pattern}, the tip forces and couples are applied along the $x_2$ and $x_3$ axes with time-varying intensities shown in Figure~\ref{fig:cerchio_FM}. The load time histories are chosen to emphasize the viscoelastic behavior of the beam.   

The total simulation time is \SI{4}{s} with a time step $dt=\SI{1E-3}{s}$. For the space discretization we used basis functions with $p=6$ and $\rm n=50$.  A material very similar to the one used in \cite{Payette&Reddy2013}, characterized by 9 Maxwell elements and $\nu=0.4$, is employed (see Table~\ref{tab:9prony}).  

\begin{table}
	\footnotesize
\centering
	\begin{tabular}{lcccc}
	 	  $E$ [$\SI{}{{N}/{m^2}}$] 	& $\tau$ [$\SI{}{s}$] \\
	\hline
\SI{1.419e9}{} &-  \\
\SI{2.977e8}{} &\SI{0.092}{}\\
\SI{6.363e7}{} &\SI{0.981}{}\\
\SI{1.583e8}{} &\SI{9.527}{}\\
\SI{1.811e8}{} &\SI{94.318}{}\\
\SI{2.388e8}{} &\SI{920.660}{}\\
\SI{2.780e8}{} &\SI{8.998e3}{}\\
\SI{3.277e8}{} &\SI{8.685e4}{}\\
\SI{3.228e8}{} &\SI{8.514e5}{}\\
\SI{4.047e8}{} &\SI{7.740e6}{}\\
	\end{tabular}
	\caption{\ang{90} circular arch: mechanical properties selected from \cite{Payette&Reddy2013}.}\label{tab:9prony}
\end{table}

\begin{figure}
\centering
	\subfigure[Tip displacement along the $x_1$-axis.\label{fig:cerchio_u1}]{\includegraphics[width=0.3\textwidth]{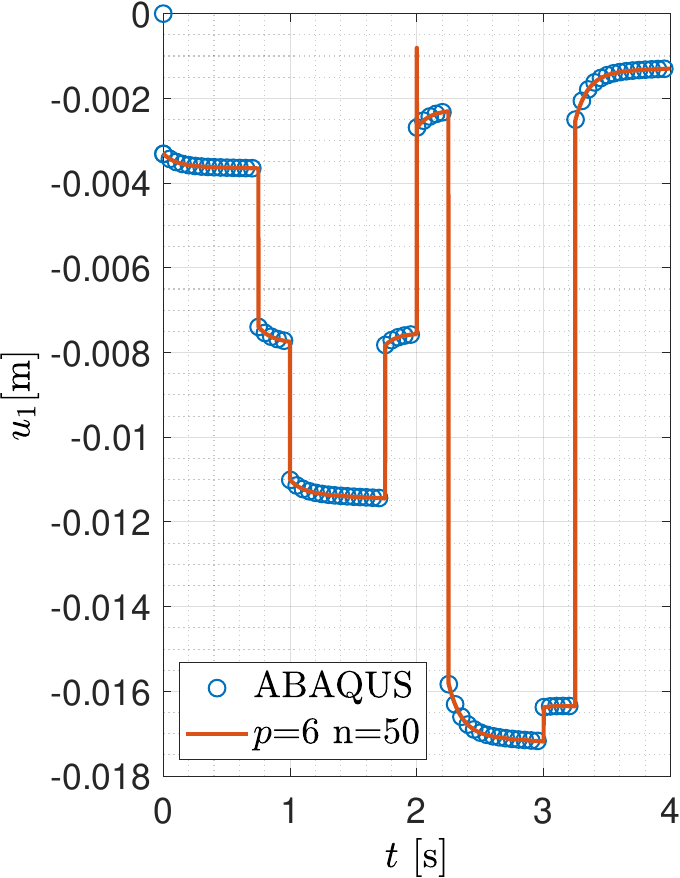}}\hspace{0.5cm}
 \subfigure[Tip displacement along the $x_2$-axis.\label{fig:cerchio_u2}]{\includegraphics[width=0.3\textwidth]{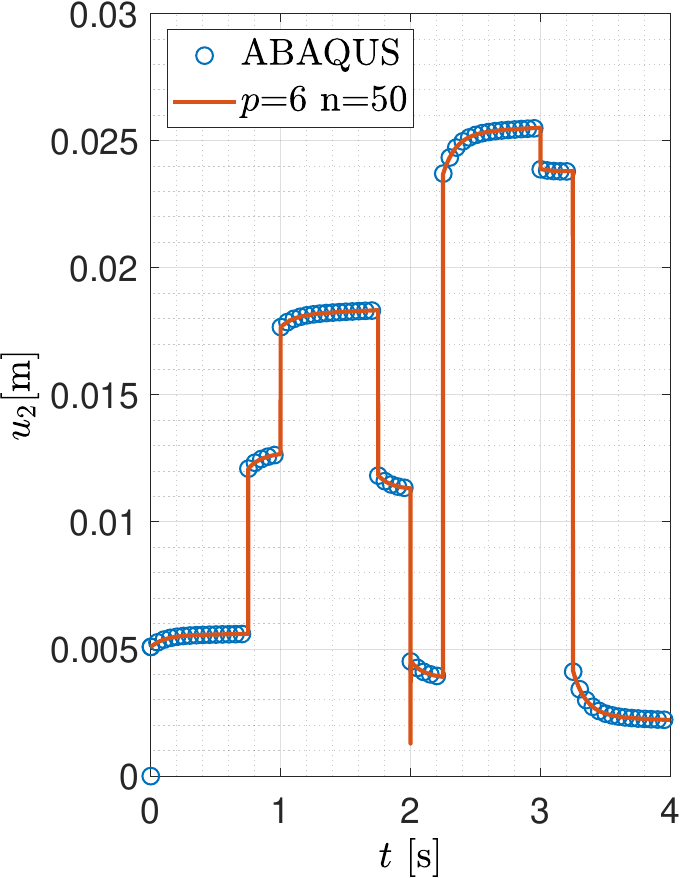}}\hspace{0.5cm}
  \subfigure[Tip displacement along the $x_3$-axis.\label{fig:cerchio_u3}]{\includegraphics[width=0.3\textwidth]{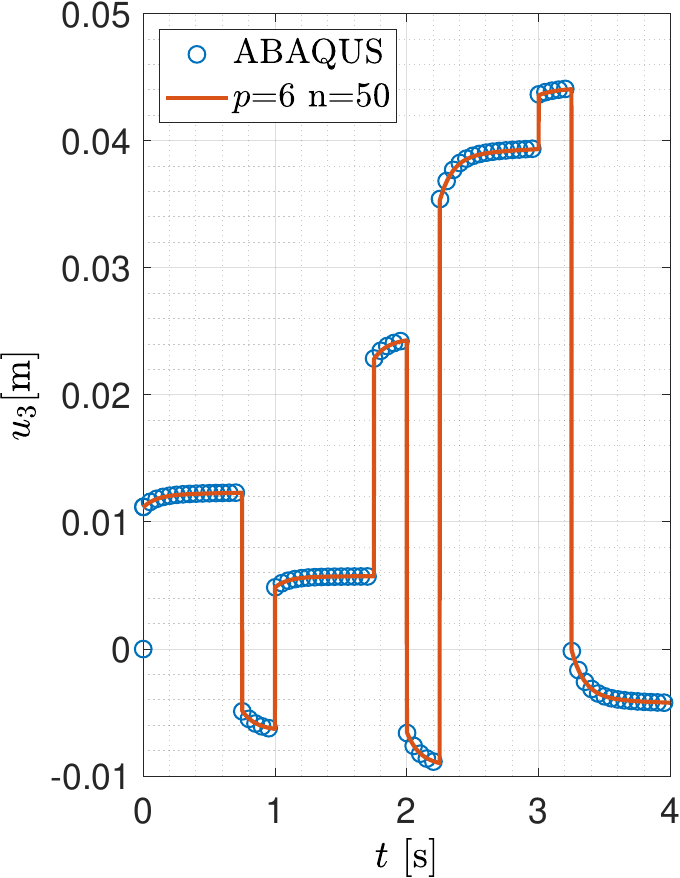}}
\caption{\label{fig:cerchio_U} \ang{90} circular arch subjected to tip forces and couples.}
\end{figure}

The tip displacements in the three Cartesian directions are shown in Figure~\ref{fig:cerchio_U}. 
A comparison with Abaqus results obtained with 314 B31 elements is also included and a very good agreement is observed. Note that the sharp variation of the displacements $u_1$ and $u_2$ at $t=\SI{2}{s}$ (see Figures~\ref{fig:cerchio_u1} and~\ref{fig:cerchio_u2}) is associated with the complete unloading of the beam for that time instant. Such a variation is not noticeable in the $x_3$ direction (see Figure~\ref{fig:cerchio_u3}) since from $t=\SI{1.999}{s}$ to $t=\SI{2.001}{s}$ $u_3$ changes sign. 

\subsection{Spivak beam under tip loads}
Moving to more complex geometrically nonlinear cases, we analyze here the Spivak beam \cite{Ignesti_etal2023}. It is a rather challenging three-dimensional geometry which permits also assessing the capability of the proposed model to reconstruct complex initial geometries featuring points with vanishing curvature. The initial beam axis is given by the following piecewise-defined curve 
\bEq
\left\{
\begin{alignedat}{3}
&\f {c}(s)=[s,\, 0,\, e^{-1/{s^2}}]\Tra \,, 	&\quad 	& s\in [-2,0)\,,\\
&\f {c}(s)=[0,\, 0,\, 0]\Tra\,,			&\quad	& s=0\,,\\
 &\f {c}(s)=[s,\, e^{-1/{s^2}},\, 0]\Tra\,,	&\quad 	& s\in (0,3].
\end{alignedat}
\right.
\eEq

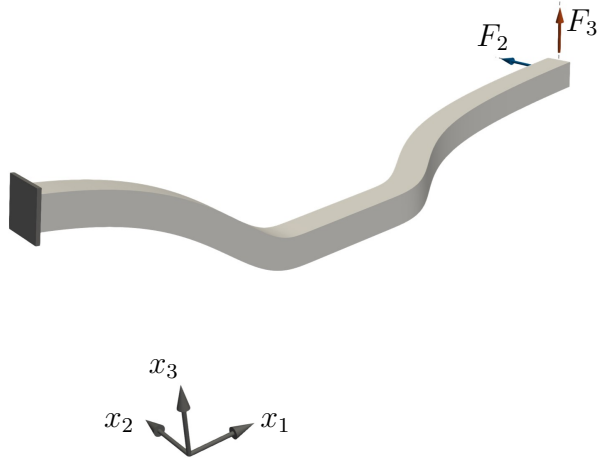
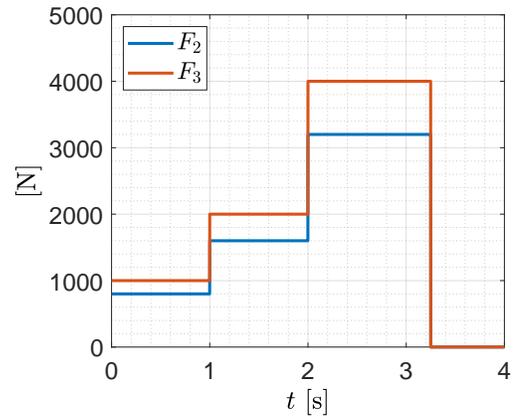
\begin{figure}
\centering
\subfigure[Beam geometry and loads.\label{fig:spivak_G}]
{\begin{overpic}[width=0.6\textwidth]
{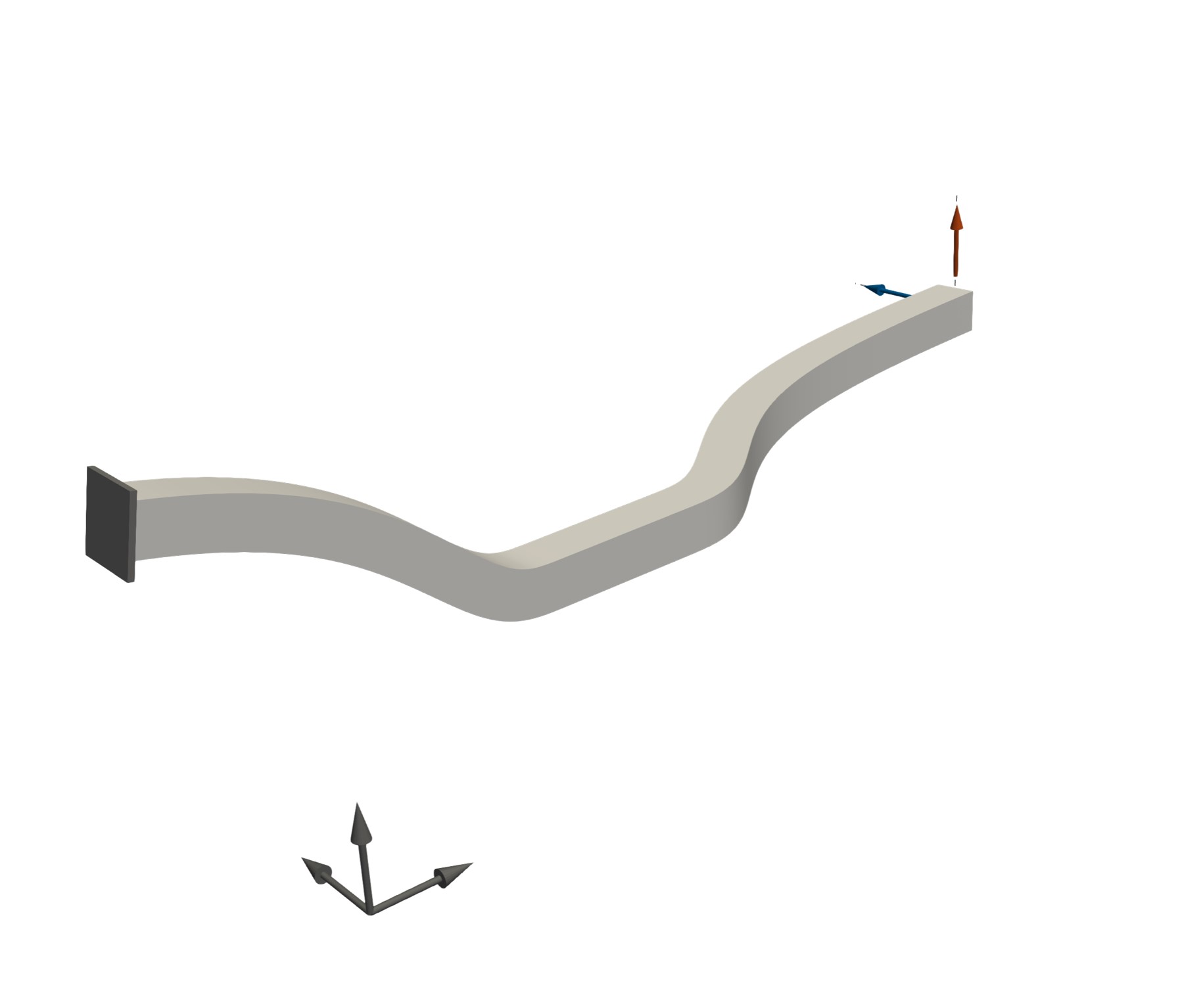}\put(410,118){$x_1$}\put(200,118){$x_2$}\put(260,190){$x_3$}
\put(700,630){$F_2$}\put(820,650){$F_3$}\end{overpic}}\hspace{-3mm}
\subfigure[Tip forces time history. \label{fig:spivak_F}]{\includegraphics[width=0.4\textwidth]{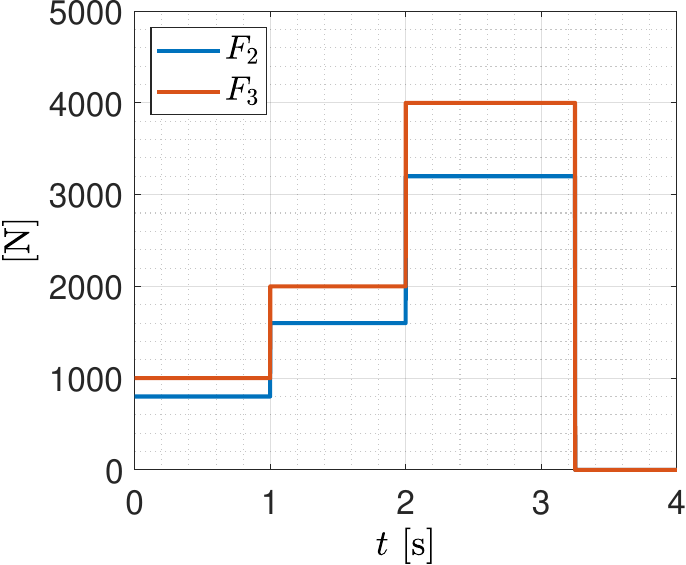}}
\caption{Spivak beam.\label{fig:spivak}}
\end{figure}

The beam has a total length of $\SI{5.50}{m}$, with a square cross section of side \SI{0.2}{m} (see Figure~\ref{fig:spivak}). The material is the same used in the previous case whose parameters are reported in Table ~\ref{tab:9prony}. The beam is loaded with two tip forces along the $x_2$ and the $x_3$ axes, respectively (see Figure~\ref{fig:spivak_G}). 
The time history of the loads is shown in Figure~\ref{fig:spivak_F}. Basis functions with $p=8$ and $\rm n=100$ are used. 

\begin{figure}
\centering
	\subfigure[Tip displacement along the $x_1$-axis.\label{fig:spivak_u1}]{\includegraphics[width=0.31\textwidth]{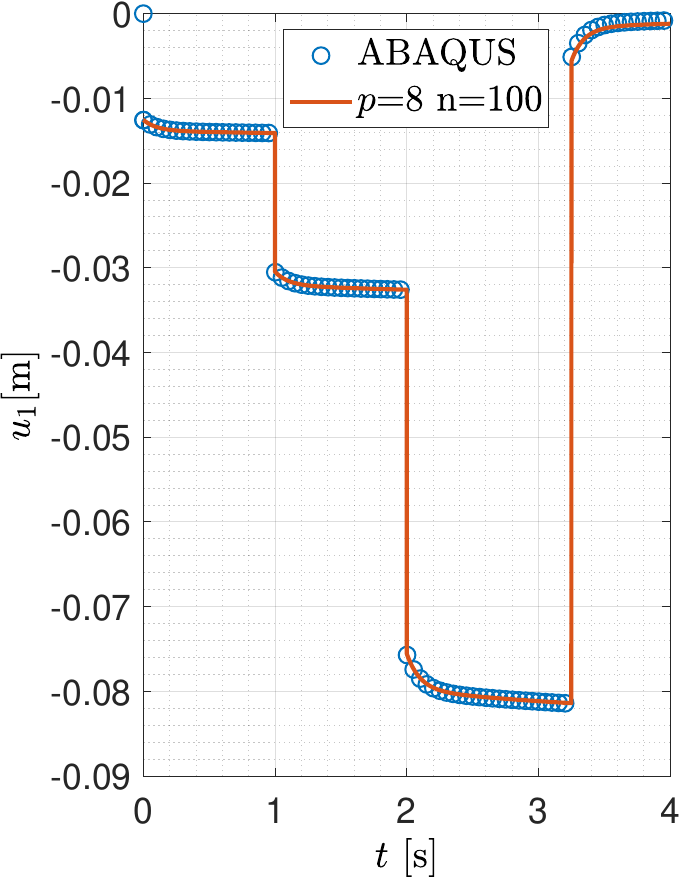}}\hspace{0.25cm}
 \subfigure[Tip displacement along the $x_2$-axis.\label{fig:spivak_u2}]{\includegraphics[width=0.31\textwidth]{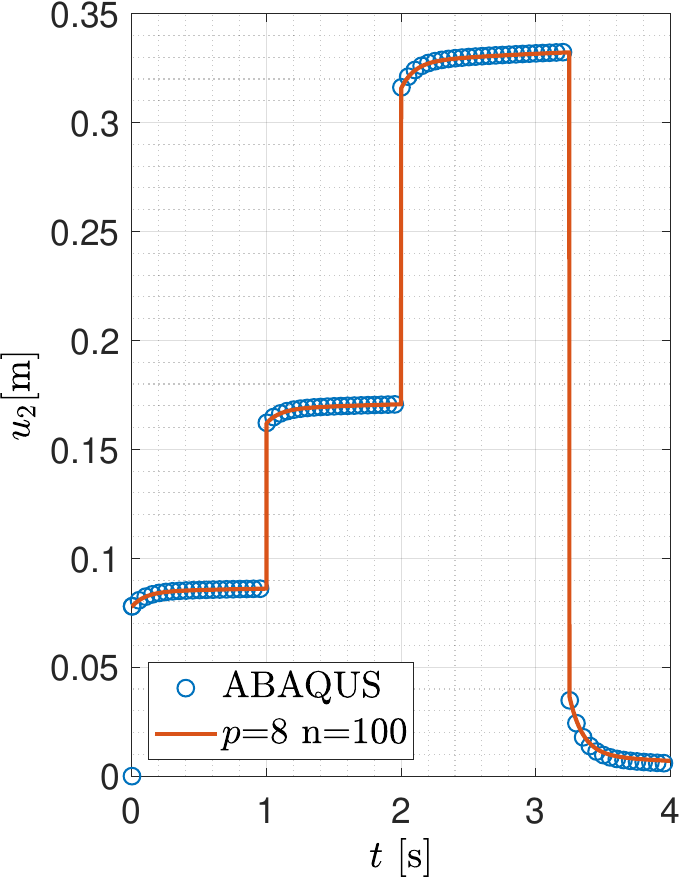}}\hspace{0.25cm}
 \subfigure[Tip displacement along the $x_3$-axis.\label{fig:spivak_u3}]{\includegraphics[width=0.31\textwidth]{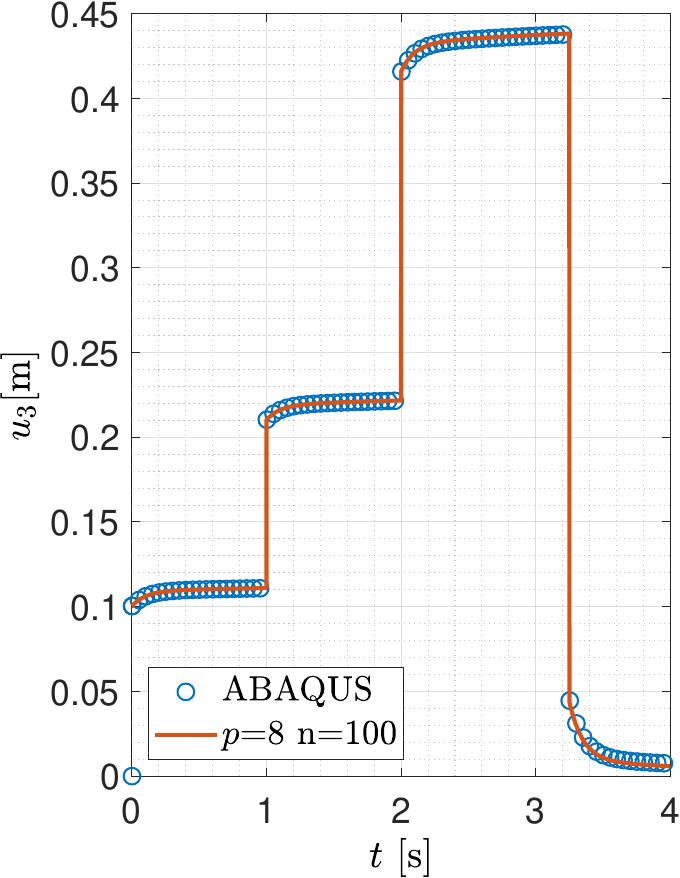}}\hspace{0.25cm}
\caption{\label{fig:spivak_U} Spivak beam results comparison.}
\end{figure}

\begin{figure}
\centering
\begin{overpic}[width=0.65\textwidth]{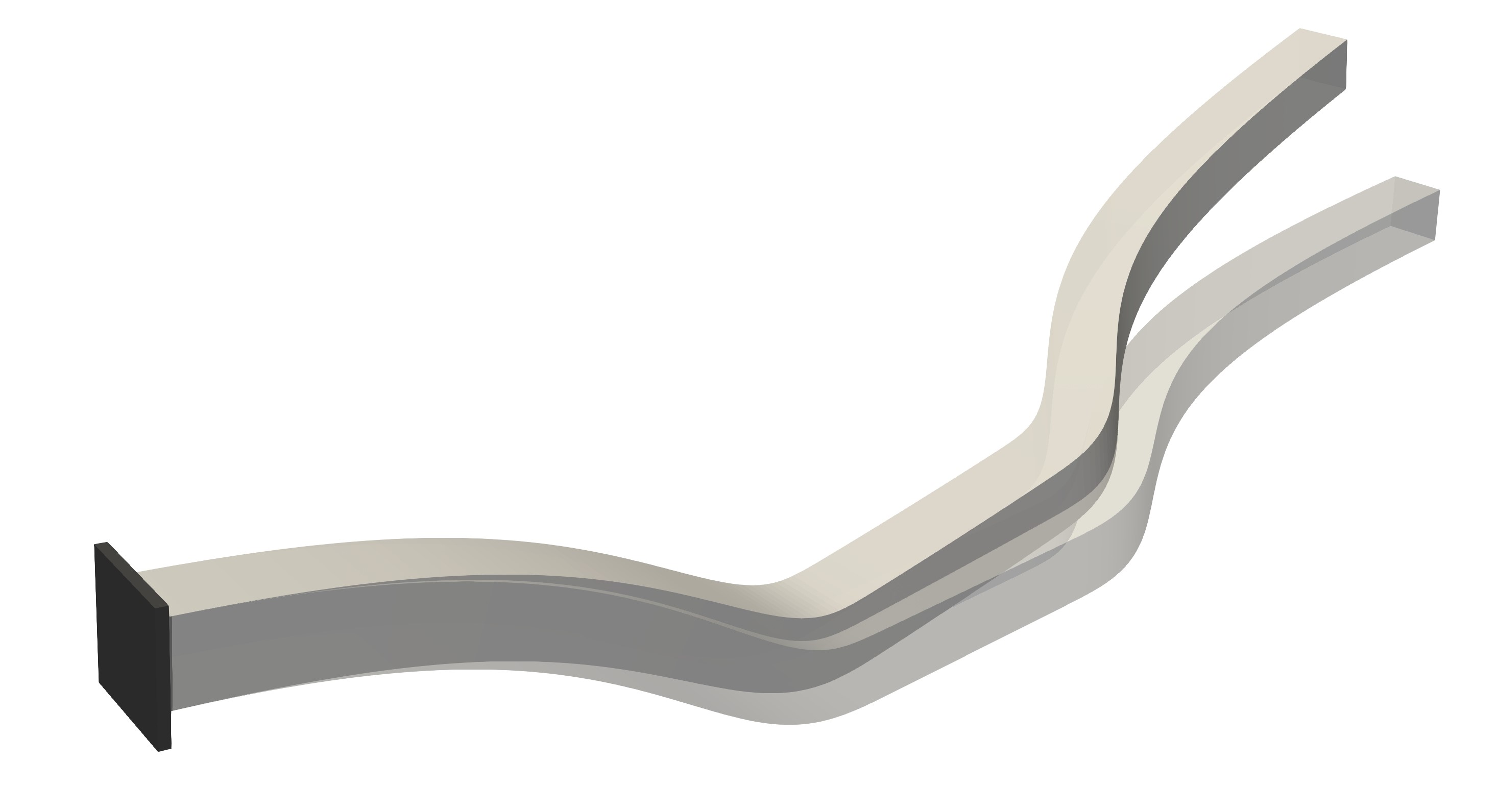}
\put(720,500){$t=\SI{3}{s}$}
\put(875,300){$t=\SI{0}{s}$}
\end{overpic}
\caption{Spivak beam clamped at one end and subjected to tip forces along $x_2$-axis and $x_3$-axis: undeformed configuration (shaded grey) deformed configuration after \SI{3}{s} (solid grey).\label{fig:spivak_def}}
\end{figure}

Also in this case, the results are compared with Abaqus, where 452 B31 finite elements are used. Figure~\ref{fig:spivak_U} shows the comparison of the tip displacements. The viscous response is very well captured during the time intervals when the loads remain constant.  
The deformed configuration at $t = \SI{3}{s}$ is shown in Figure~\ref{fig:spivak_def}. 

\subsection{Lissajous beam}
\begin{figure}
\begin{minipage}{0.65\textwidth}
\centering
\subfigure[Beam geometry and loads.\label{fig:Lissa_G}]{\begin{overpic}[clip,width=0.8\textwidth]{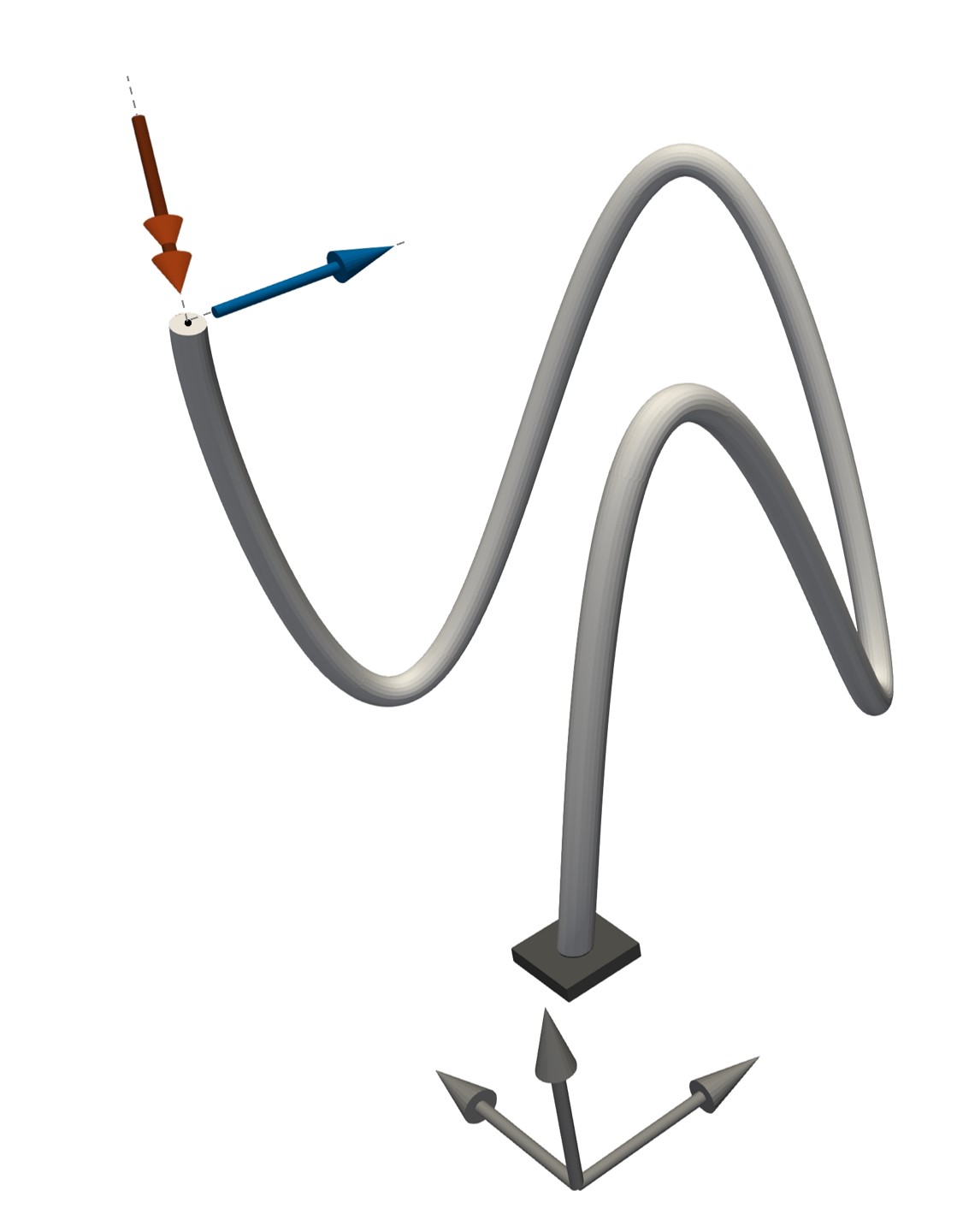}\put(625,125){$x_1$}\put(300,125){$x_2$}\put(390,150){$x_3$}\put(200,800){$F_1$}\put(40,850){$M_3$}\end{overpic}}
\end{minipage}
\begin{minipage}{0.35\textwidth}
\centering
\subfigure[Time history of the tip force, $F_1$.\label{fig:Lissa_F}]{\includegraphics[width=1\textwidth]{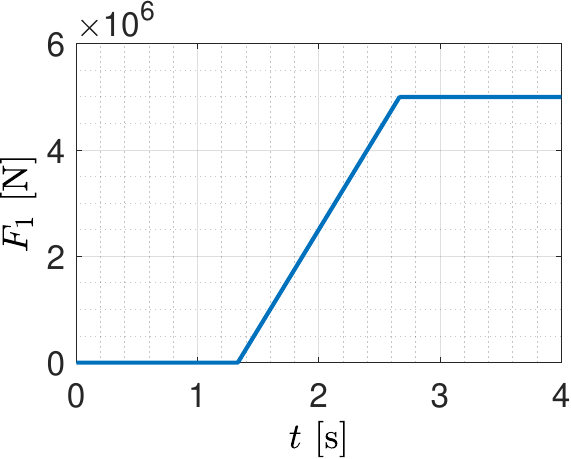}}
\subfigure[Time history of the tip couple, $M_3$.\label{fig:Lissa_M}]{\includegraphics[width=1\textwidth]{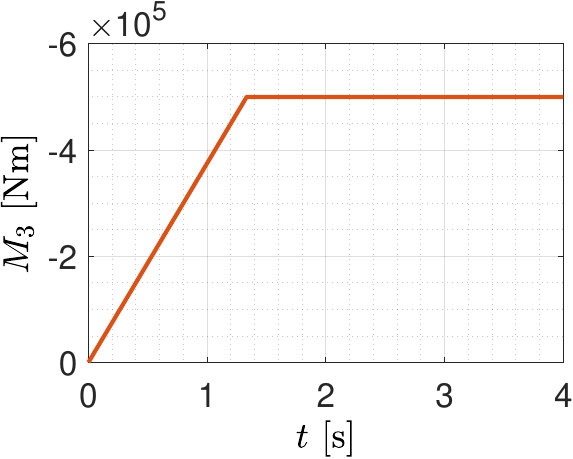}}
\end{minipage}
\caption{Lissajous beam.\label{fig:Lissa}}
\end{figure}

The Lissajous beam is another geometrically challenging case useful to test the capabilities to model viscoelastic beams with repeated strong curvatures variations. The initial beam axis is defined by the curve 
$\f c(s) = [\cos(3s),\, \sin(2s),\, \sin(7s)]\Tra$ with $s\in [-\pi/3, \pi/3]$ (see Figure~\ref{fig:Lissa}). 
The beam has a circular cross section of diameter \SI{0.12}{m}. 
The viscoelastic material is represented by one Maxwell element with parameters shown in Table~\ref{tab:PLA}. A constant Poisson ratio $\nu=0.4$ is assumed. These material properties are selected from those of the Polylactic acid (PLA) studied in \cite{Wan_etal2022}.

The beam is clamped at one end and is loaded at its tip with a negative couple in the $x_3$ direction and a positive force in the $x_1$ direction (see Figures~\ref{fig:Lissa_F} and \ref{fig:Lissa_M}). 

The test is carried with $\rm n=100$ and $p=8$. The total simulation time is $\SI{4}{s}$ with a time step $dt=\SI{1E-3}{s}$.

\begin{table}
	\footnotesize
\centering
	\begin{tabular}{lcccc}
$E$ [$\SI{}{{N}/{m^2}}$] 	& $\tau$ [$\SI{}{s}$] \\
	\hline
\SI{2.80e5}{} &-   \\
\SI{3.61e7}{} &\SI{0.2}{}\\
	\end{tabular}
		\caption{Mechanical properties of PLA modelled with a single Maxwell element \cite{Wan_etal2022}.}\label{tab:PLA}
\end{table}

\begin{figure}
\centering
\subfigure[Deformation at $t=\SI{1.33}{s}$: perspective view.\label{fig:lissa_133d}]{\begin{overpic}[width=0.49\textwidth]{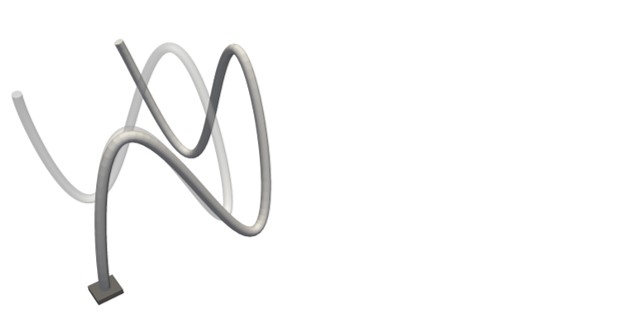}
\end{overpic}}
\subfigure[Deformation at $t=\SI{1.33}{s}$: top view.\label{fig:lissa_133t}]{\begin{overpic}[width=0.49\textwidth]{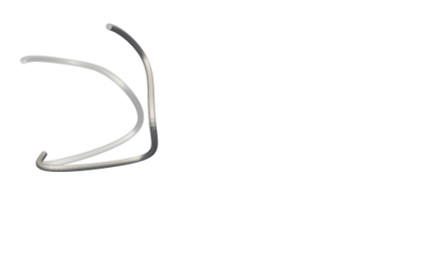}
\end{overpic}}
\subfigure[Deformation at $t=\SI{2.67}{s}$: perspective view.\label{fig:lissa_267d}]{\begin{overpic}[width=0.49\textwidth]{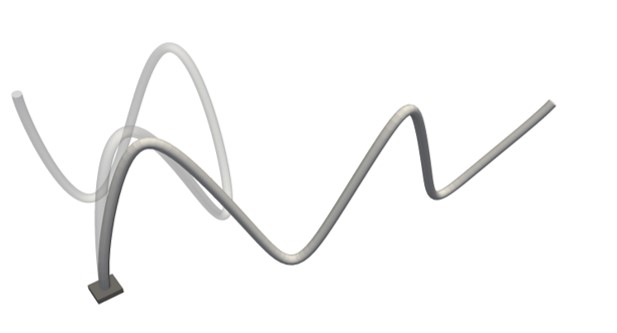}
\end{overpic}}
\subfigure[Deformation at $t=\SI{2.67}{s}$: top view.\label{fig:lissa_267t}]{\begin{overpic}[width=0.49\textwidth]{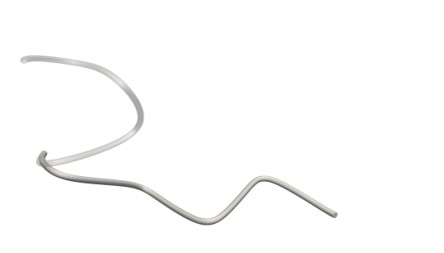}
\end{overpic}}
\subfigure[Deformation at $t=\SI{4.00}{s}$: perspective view.\label{fig:lissa_400d}]{\begin{overpic}[width=0.49\textwidth]{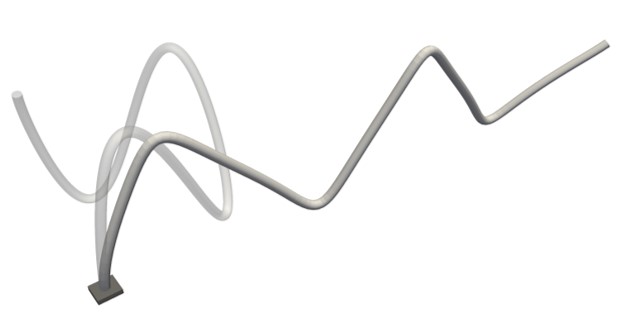}
\end{overpic}}
\subfigure[Deformation at $t=\SI{4.00}{s}$: top view.\label{fig:lissa_400t}]{\begin{overpic}[width=0.49\textwidth]{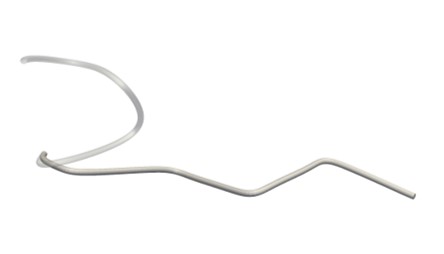}
\end{overpic}}
\caption{Selected snapshots of the complex deformation of the Lissajous beam.\label{fig:lissa_def}}
\end{figure}

Results are reported in Figure~\ref{fig:lissa_def} (left column: perspective views, right column: top views). Figures~\ref{fig:lissa_133d} and \ref{fig:lissa_133t} show the deformed beam configuration at the end of the ramp of the concentrated moment $M_3$ ($t=\SI{1.33}{s}$). The beam tends to twist around $x_3$-axis with respect to its clamped end. Figures~\ref{fig:lissa_267d} and \ref{fig:lissa_267t} show the beam configuration at $t=\SI{2.67}{s}$, where it is noticeable the joint effect of $F_1$ and $M_3$ in straightening the beam. 
Moreover, it can also be noticed the progressive twisting of the beam around $x_3$-axis with respect to Figures~\ref{fig:lissa_133d} and \ref{fig:lissa_133t}. The final configuration is reported in Figures~\ref{fig:lissa_400d} and \ref{fig:lissa_400t}. From $t=\SI{2.67}{s}$ to $t=\SI{4}{s}$, both 
$F_1$ and $M_3$ are kept constant (see Figure~\ref{fig:Lissa_F} and~\ref{fig:Lissa_M}), therefore the additional displacements, observable by comparing Figure~\ref{fig:lissa_400d} with Figure~\ref{fig:lissa_267d}, are completely ascribed to viscous effects.

\subsection{Planar and tubolar nets with curved cell elements}
\begin{figure}
\centering
\subfigure[Single cell shape, 3D view and loadings.\label{fig:13a}]{\begin{overpic}[width=0.75\textwidth]{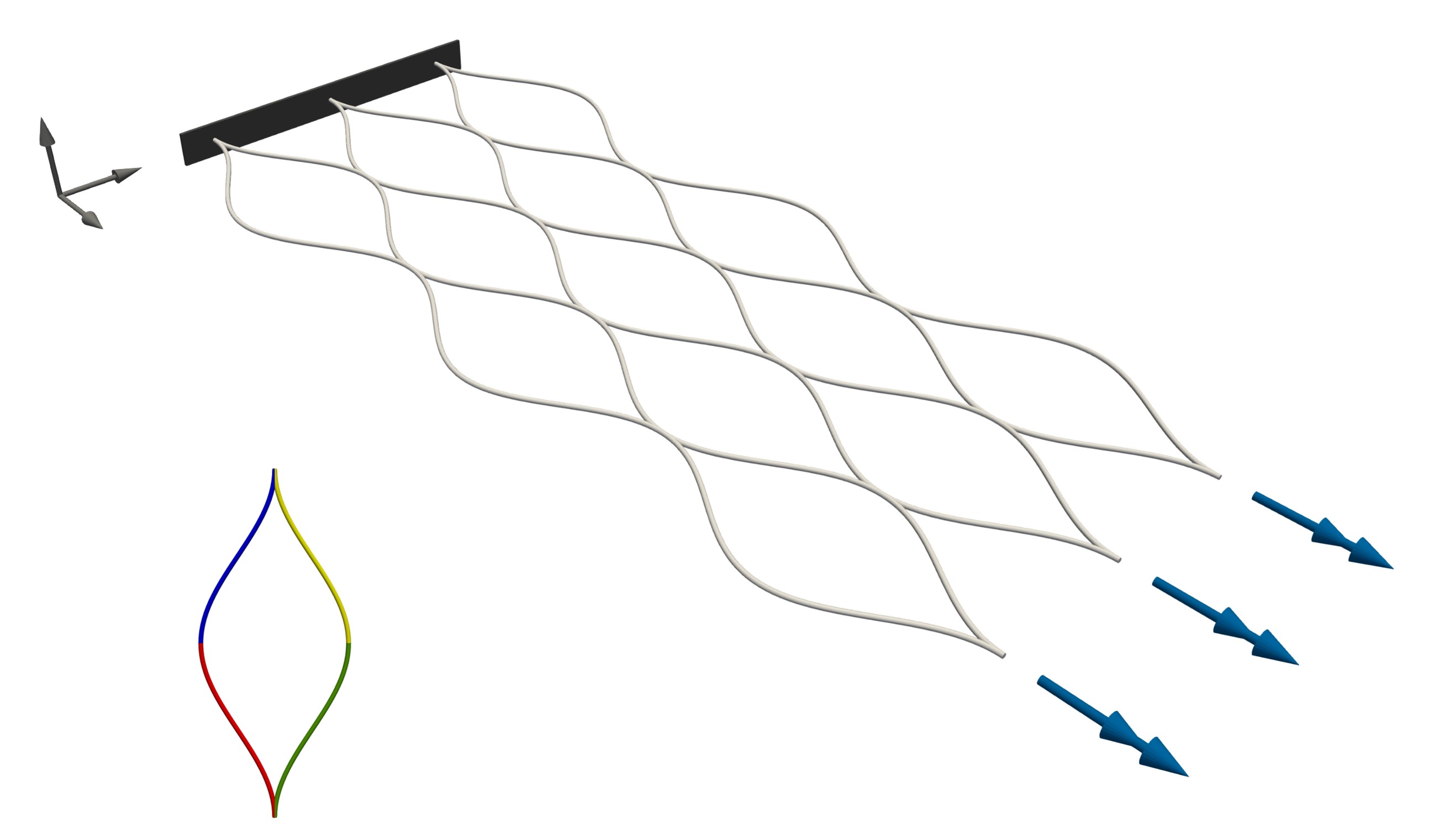}\put(25,400){$x_1$}\put(70,480){$x_2$}\put(0,510){$x_3$}\put(790,85){$M_1$}\put(870,160){$M_1$}\put(935,225){$M_1$}\put(235,60){I}\put(230,190){II}\put(100,60){IV}\put(100,190){III}\end{overpic}}\hspace{0.25cm}
\subfigure[Plane view of the net (units in \SI{}{cm}).\label{fig:13b}]{\begin{overpic}[width=0.20\textwidth]{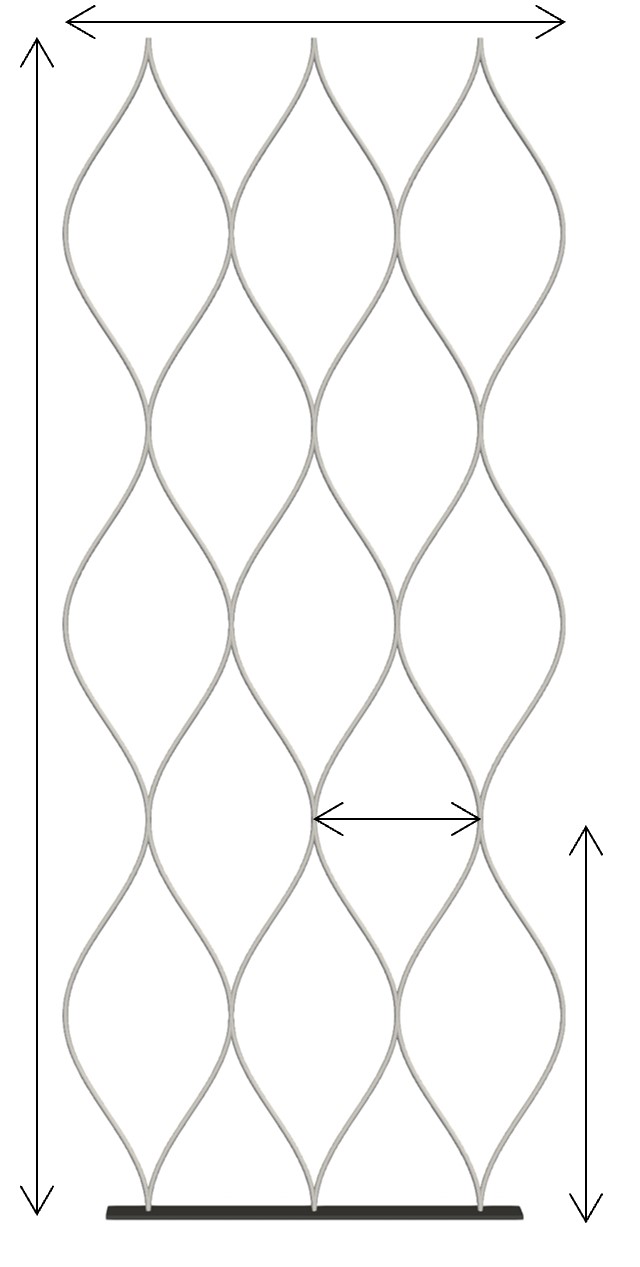}\put(285,350){\colorbox{white}{\makebox(20,20){\textcolor{black}{4}}}}\put(475,150){$3\pi$\colorbox{white}}\put(0,350){\colorbox{white}{\makebox(20,20){\textcolor{black}{$9\pi$}}}}\put(220,1000){\colorbox{white}{\makebox(20,20){\textcolor{black}
{$12$}}}}\end{overpic}}
\subfigure[Time-varying loads applied to the free tips of the net.\label{fig:13c}]{\includegraphics[width=0.5\textwidth]{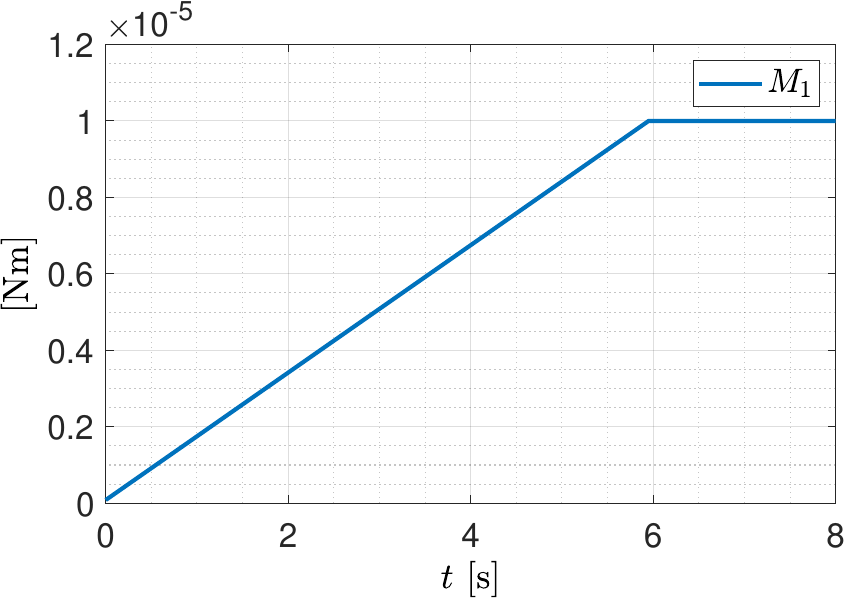}}
\caption{Twisting of a planar net made of curved cell elements.\label{fig:13}}
\end{figure}

\begin{figure}
\centering
\subfigure[3D view.\label{fig:14a}]{\begin{overpic}[width=0.47\textwidth]{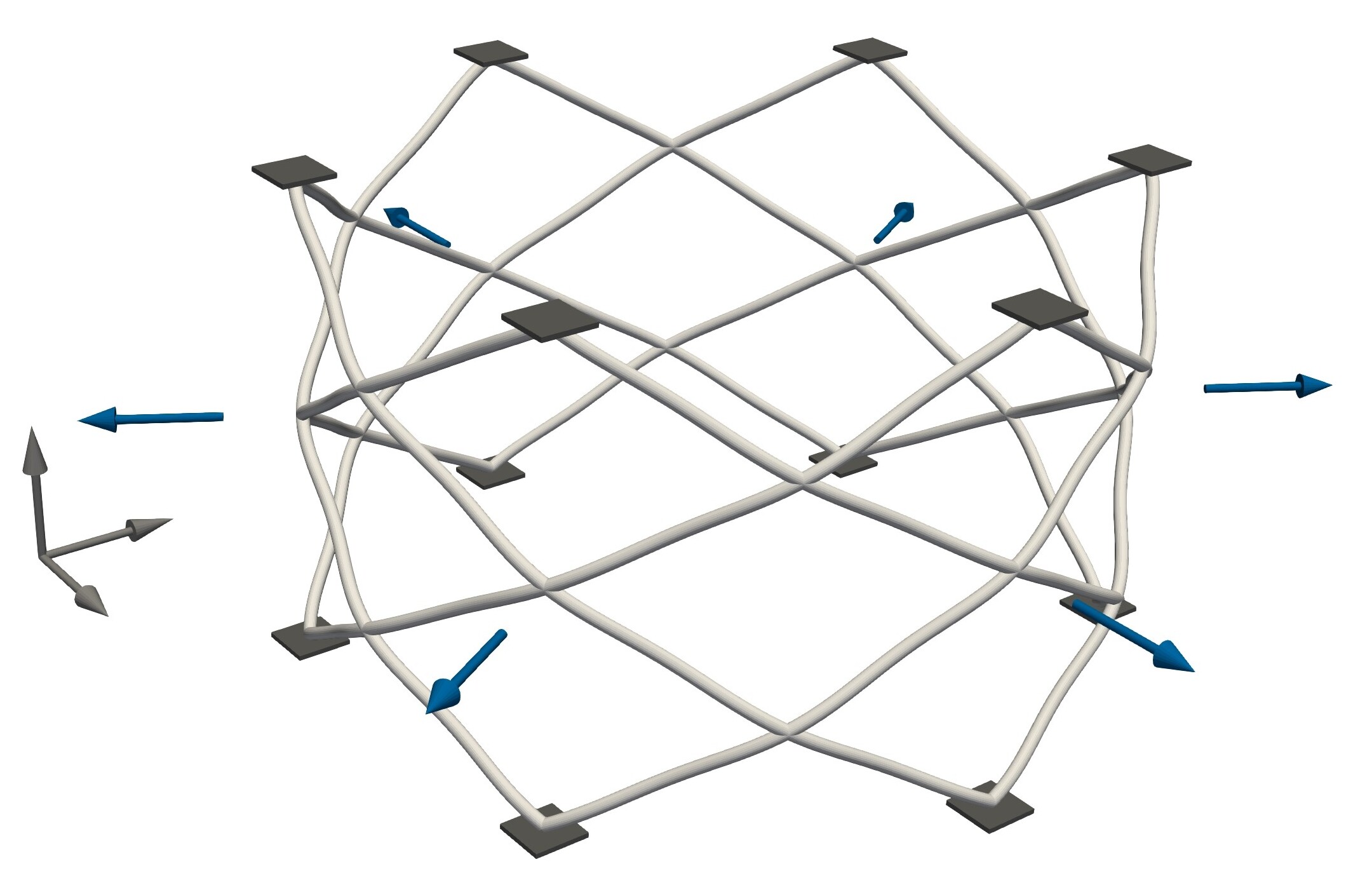}\put(25,160){$x_1$}\put(100,270){$x_2$}\put(0,345){$x_3$}
\end{overpic}}\hspace{0.25cm}
\subfigure[Lateral view in the $x_2$-$x_3$ plane.\label{fig:14b}]{\begin{overpic}[width=0.47\textwidth]{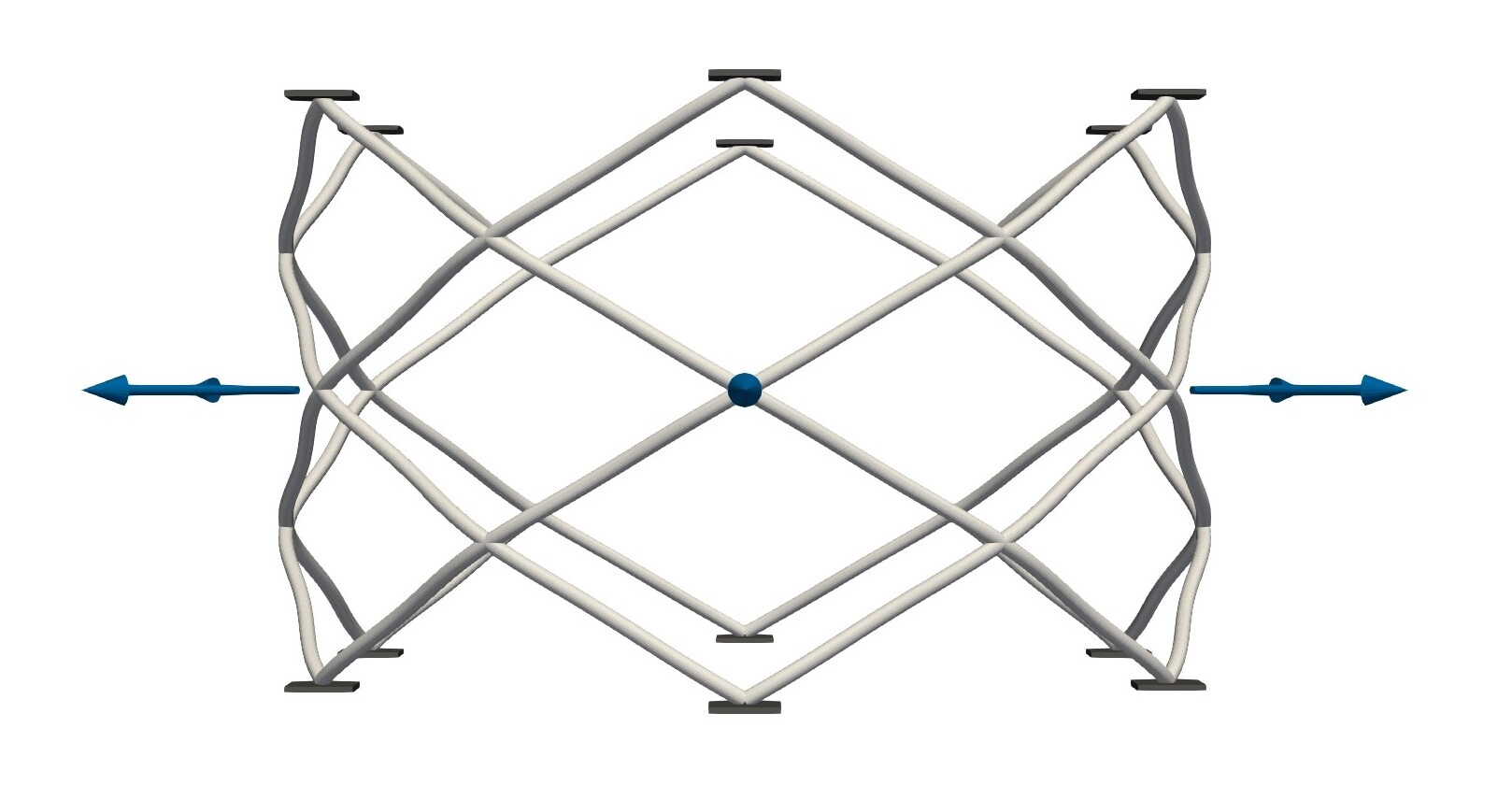}
\end{overpic}}\hspace{0.25cm}
\subfigure[Lateral view in the $x_1$-$x_3$ plane.\label{fig:14c}]{\begin{overpic}[width=0.47\textwidth]{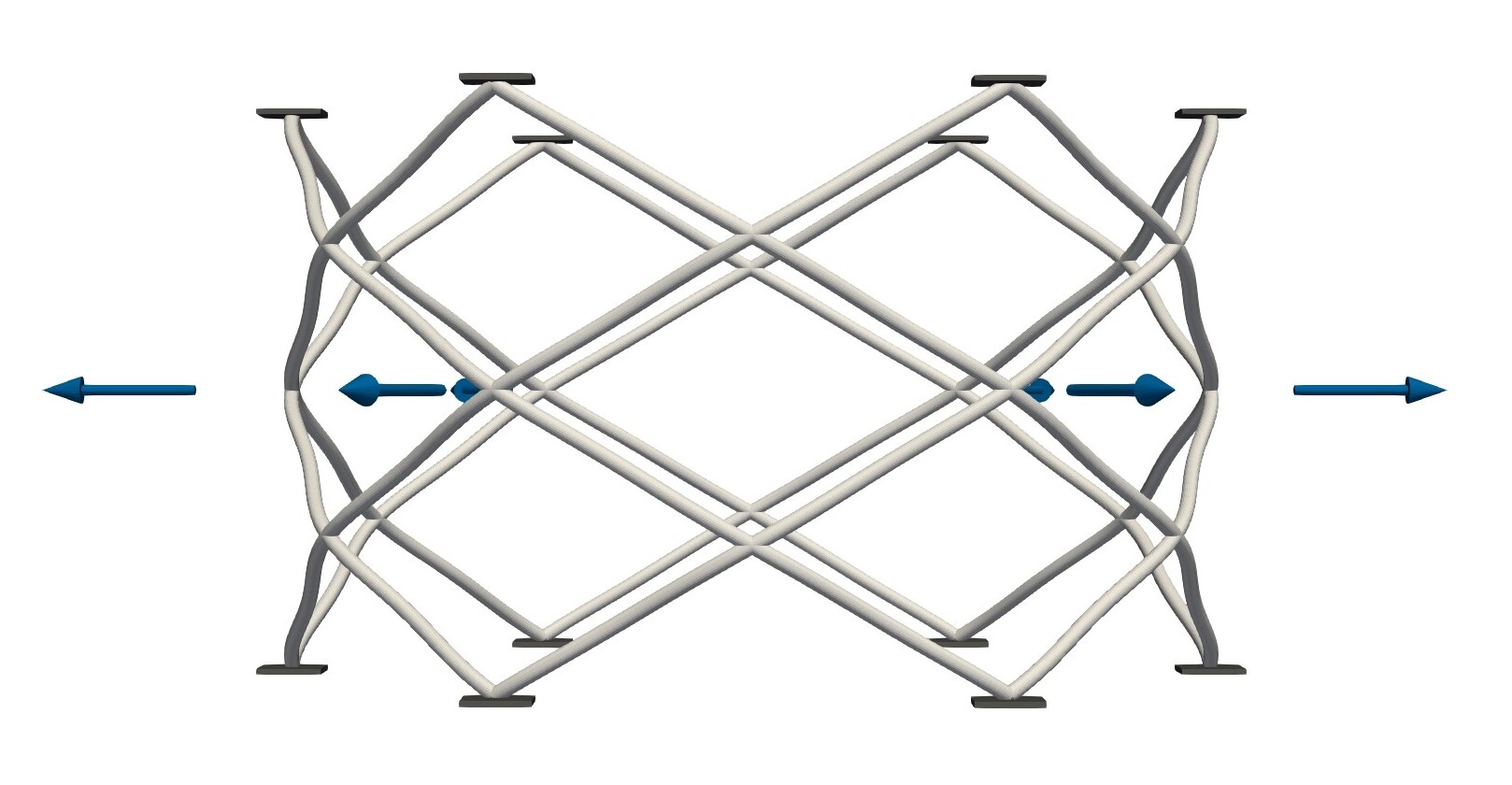}
\end{overpic}}\hspace{0.25cm}
\subfigure[Top view in the $x_1$-$x_2$ plane.\label{fig:14d}]{\begin{overpic}[width=0.47\textwidth]{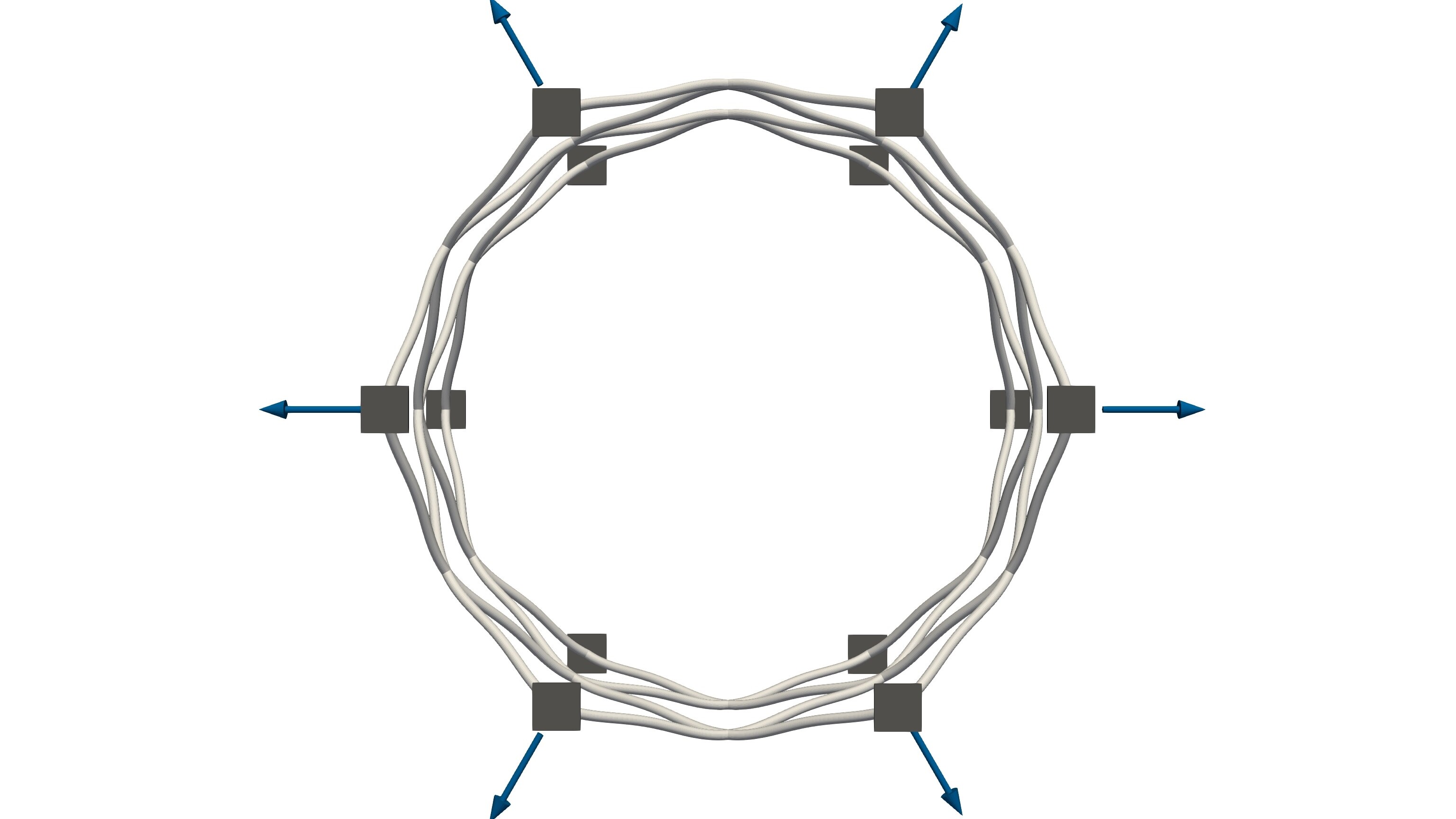}
\end{overpic}}\hspace{0.25cm}
\caption{Tubular net subjected to concentrated radial forces.\label{fig:14}}
\end{figure}

In the last numerical tests, we address the cases of complex multi-patch structures whose cells are made of curved beam elements. We first consider the case of a planar net and then a cylindrical one (see Figures~\ref{fig:13} and \ref{fig:14}). For both cases, the same material of the Lissajous beam is employed with a circular cross section of radius \SI{0.65}{mm}.

The planar net (Figure~\ref{fig:13}) is composed by 36 patches forming a network of 13 closed curved cells (see bottom left of Figure~\ref{fig:13a}). The centroid line of each patch is a curve described by $\f c(s) = [1.5s,\, -1^{k}\sin(s-\pi/2)+2(k-1),\, 0]\Tra$, where $s\in [(j-1)\pi, j\pi]$, $k=1,2,...,n_1$ represents the number of subdivisions in the $x_1$-direction and $j=1,2,...,n_2$ the subdivisions in the $x_2$-direction. 
Each patch is discretized with basis functions with $p=8$ and $\rm n=35$.
The structure is clamped at one end it is loaded with concentrated couples along the $x_1$-axis along the opposite side. These couples collectively increase linearly from $\SI{0}{s}$ to $\SI{6}{s}$, then they remain constant until $\SI{8}{s}$ (see Figure~\ref{fig:13c}). 

The second structure is a tubular net (see Figure~\ref{fig:14}) consisting in 48 patches connected together to form 12 helical wires. Each patch is described, in a polar coordinate system with $\vtht \in [0,\pi]$, by the following curve
\bEq
\f c(s) = [(R+r\cos(\vtht n_w)) \cos(\vtht) ,\, (R+r\cos(\vtht n_w))\sin(\vtht)  ,\, R\vtht\tan(\pi/6)]\Tra\,,
\eEq
where $n_w=4$ denotes the number of patch per wire, while $R=\SI{4e-2}{m}$ and $r=\SI{6.5e-4}{m}$ are the outer radius of the tube and the radius of the wires, respectively.
Each patch is discretized with basis functions with $p=8$ and $\rm n=20$. The beams system is clamped at both ends and loaded by concentrated radial forces at the nodes of the central section (see Figure~\ref{fig:14}). The loads are applied with a linear ramp for $\SI{1.15}{s}$, and then kept constant at $\SI{1}{N}$ until the end of the simulation. The simulation time is $\SI{6}{s}$, with $dt=\SI{5e-2}{s}$. 

\begin{figure}
\centering
\subfigure[$t=\SI{2}{s}$.\label{fig:16a}]{\begin{overpic}[clip,width=0.49\textwidth]{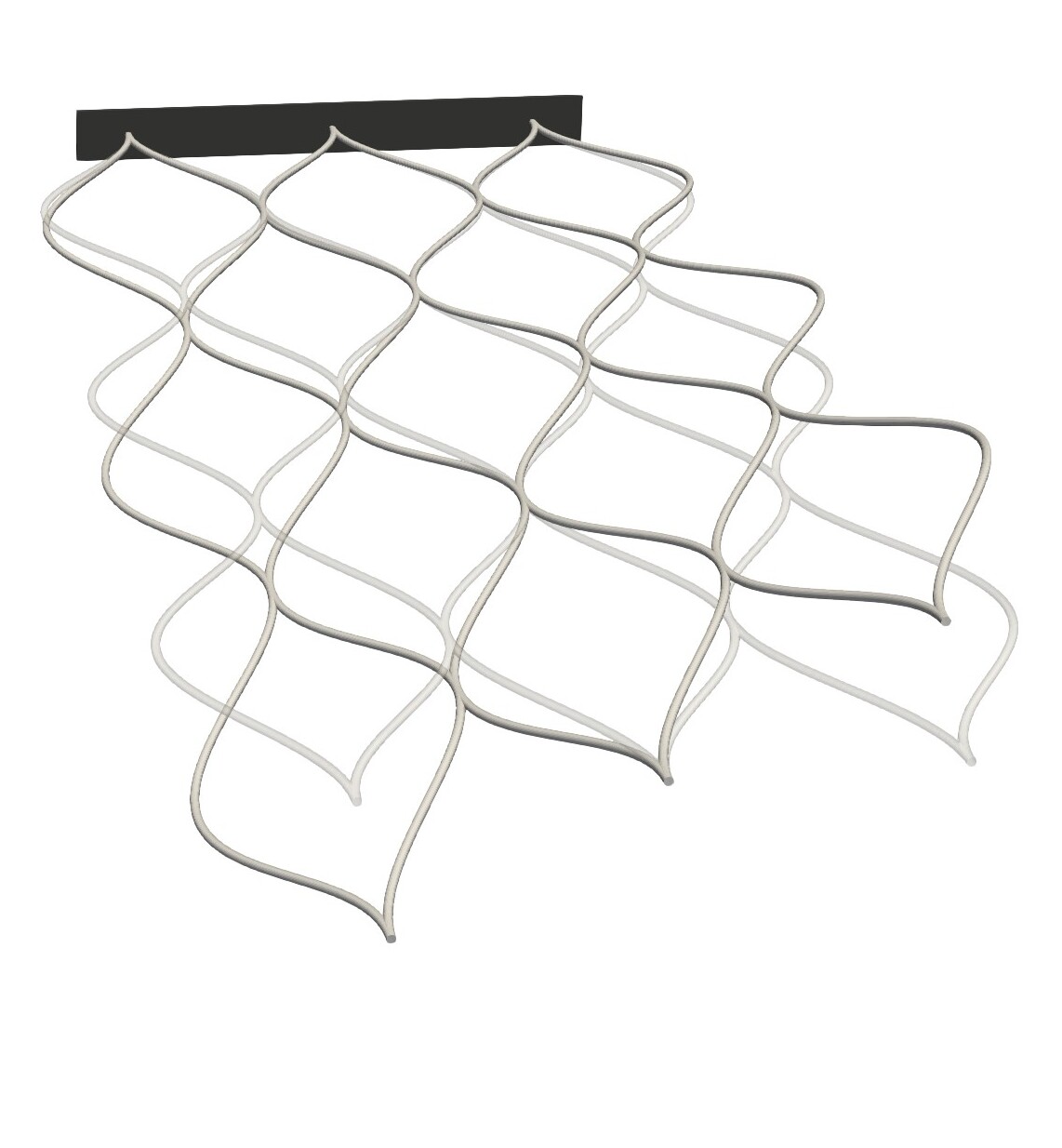}
\end{overpic}}
\subfigure[$t=\SI{4}{s}$.\label{fig:16a}]{\begin{overpic}[clip,width=0.49\textwidth]{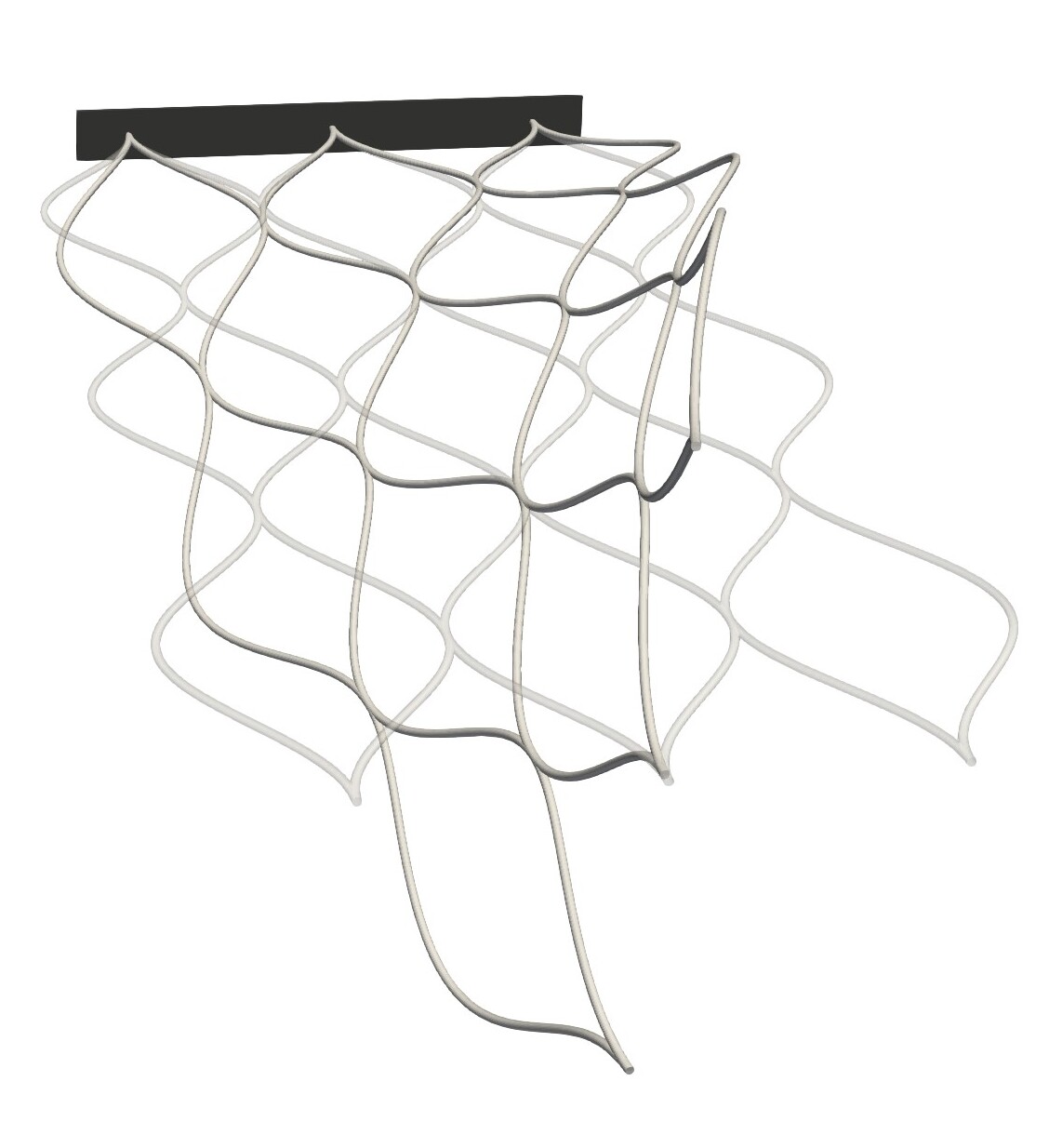}
\end{overpic}}
\subfigure[$t=\SI{6}{s}$.\label{fig:16a}]{\begin{overpic}[width=0.49\textwidth]{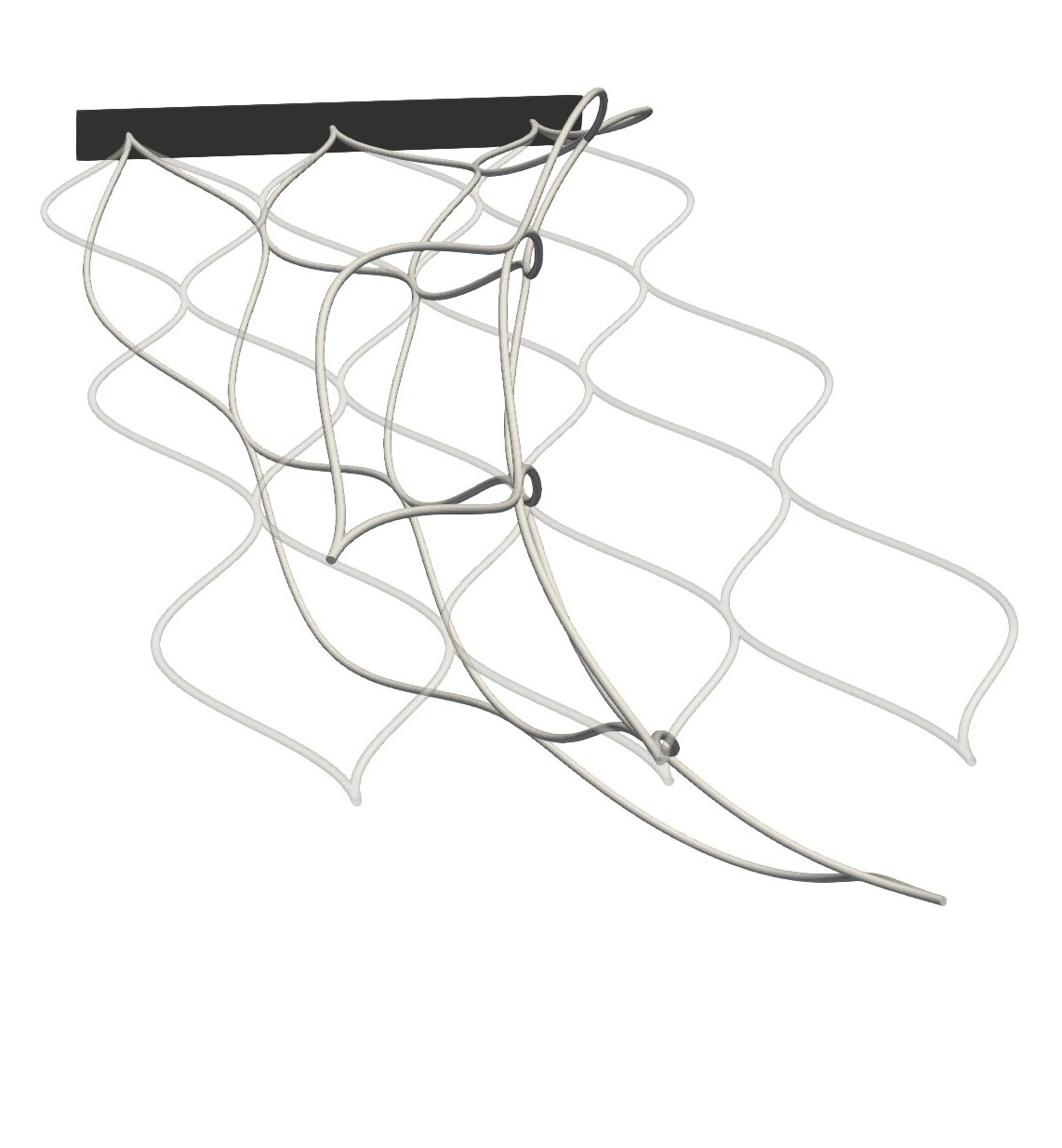}
\end{overpic}}
\subfigure[$t=\SI{8}{s}$.\label{fig:16a}]{\begin{overpic}[width=0.49\textwidth]{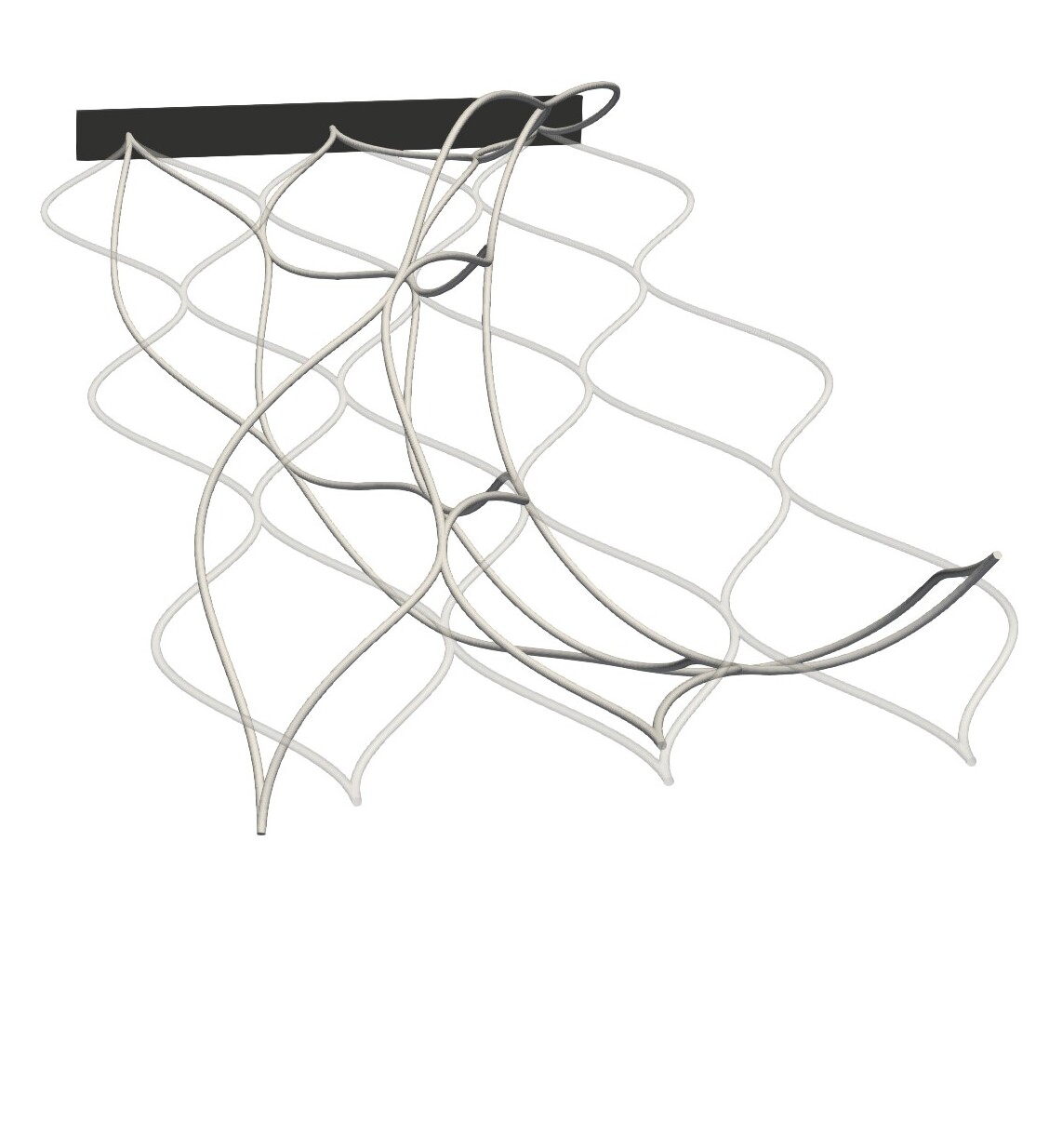}
\end{overpic}}
\caption{Planar net: 3D views of the deformed shape at different time instants.\label{fig:15}}
\end{figure}

\begin{figure}
\centering
\subfigure[Top and front view, $t=\SI{2}{s}$.\label{fig:16a}]{\begin{overpic}[width=1\textwidth]{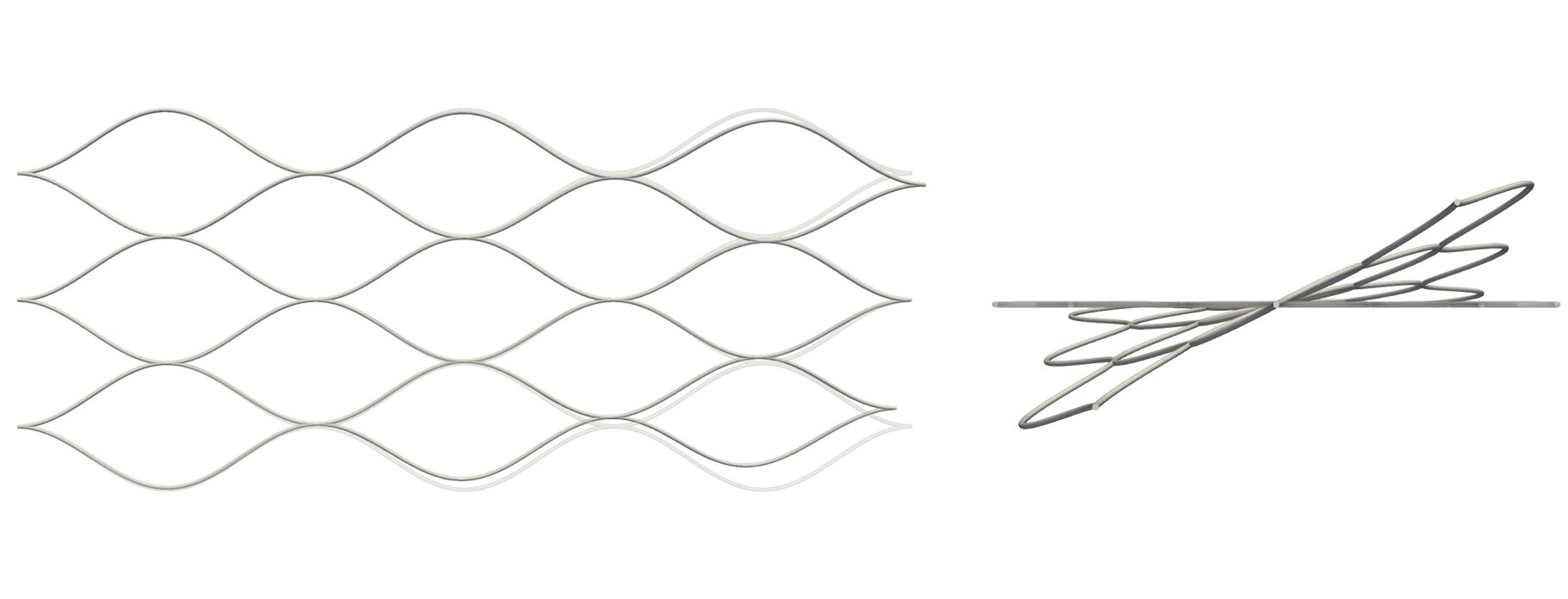}
\end{overpic}}
\subfigure[Top and front view, $t=\SI{4}{s}$.\label{fig:16b}]{\begin{overpic}[width=1\textwidth]{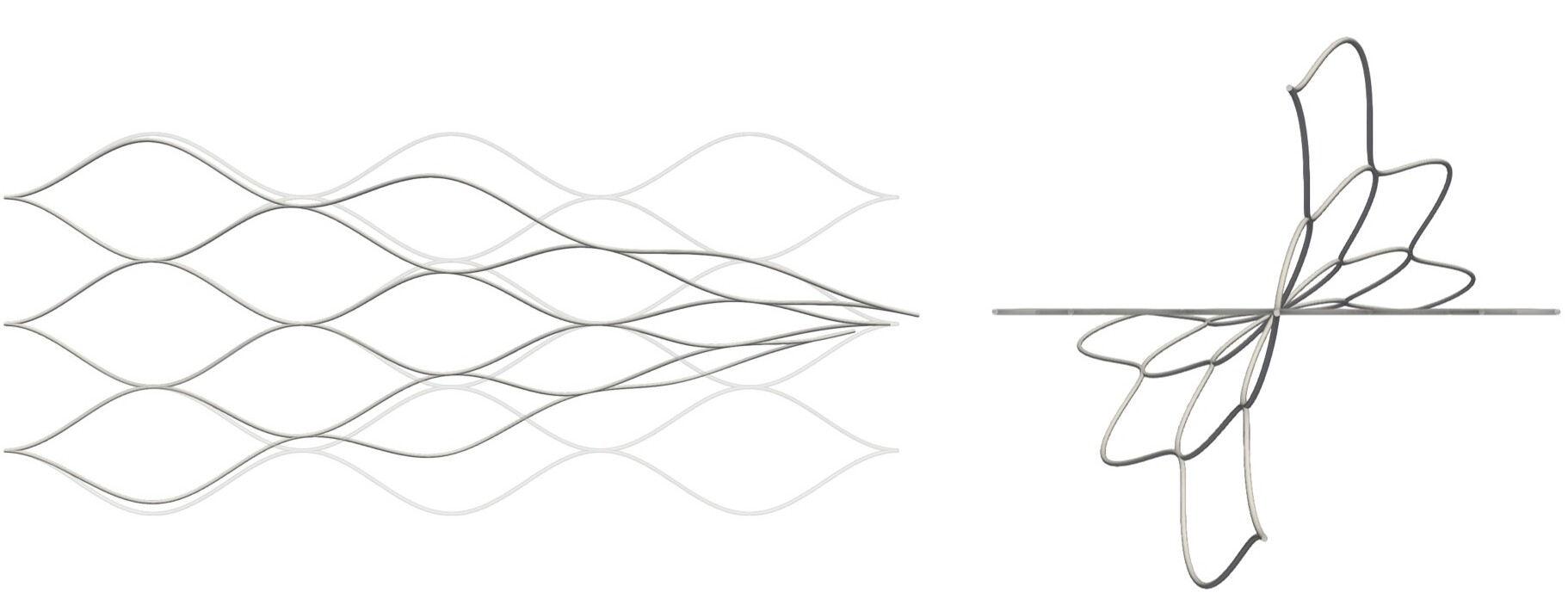}
\end{overpic}}
\subfigure[Top and front view, $t=\SI{6}{s}$.\label{fig:16c}]{\begin{overpic}[width=1\textwidth]{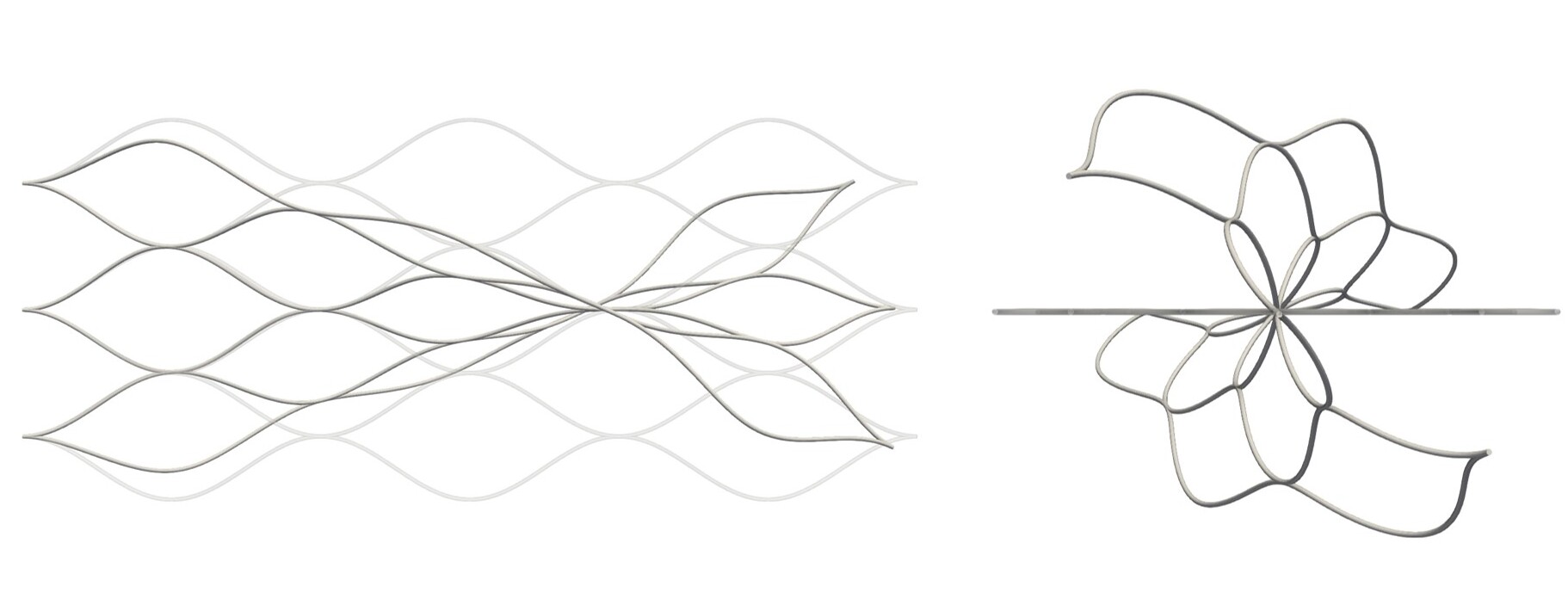}
\end{overpic}}
\end{figure}

\begin{figure}
\subfigure[Top and front view, $t=\SI{8}{s}$.\label{fig:16d}]{\begin{overpic}[width=1\textwidth]{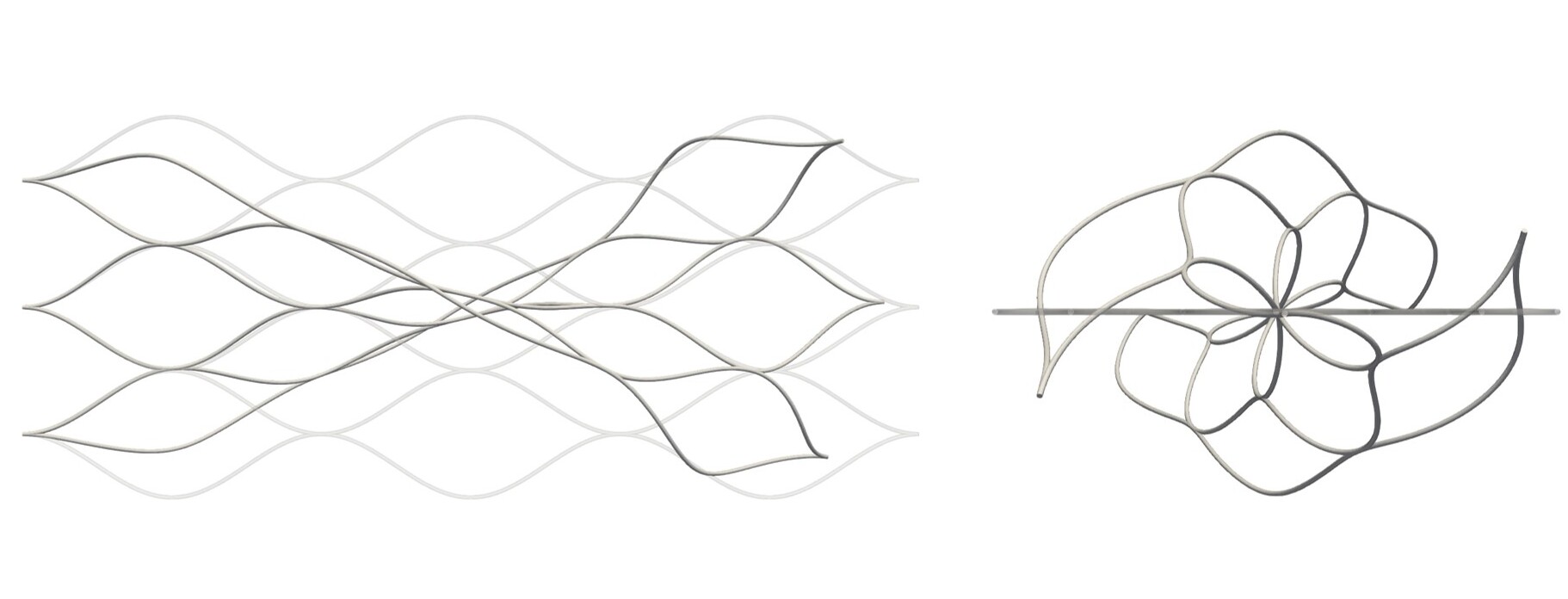}
\end{overpic}}
\caption{Planar net: plane views of the deformed shape at different time instants.\label{fig:16}}
\end{figure}

Results for the planar net are reported in Figures~\ref{fig:15} and \ref{fig:16} for $t=\SI{2}{s}$, \SI{4}{s}, \SI{6}{s}, and \SI{8}{s}. Significant viscous deformations are noticeable from $t=\SI{6}{s}$ to $t=\SI{8}{s}$, that lead to the complete twist of the net with rotations larger than $\pi$ (see Figures~\ref{fig:16c} and \ref{fig:16d}).  

\begin{figure}
\centering
\subfigure[View on $x_2$-$x_3$, $t=\SI{1.15}{s}$.\label{fig:17a}]{\begin{overpic}[width=0.32\textwidth]{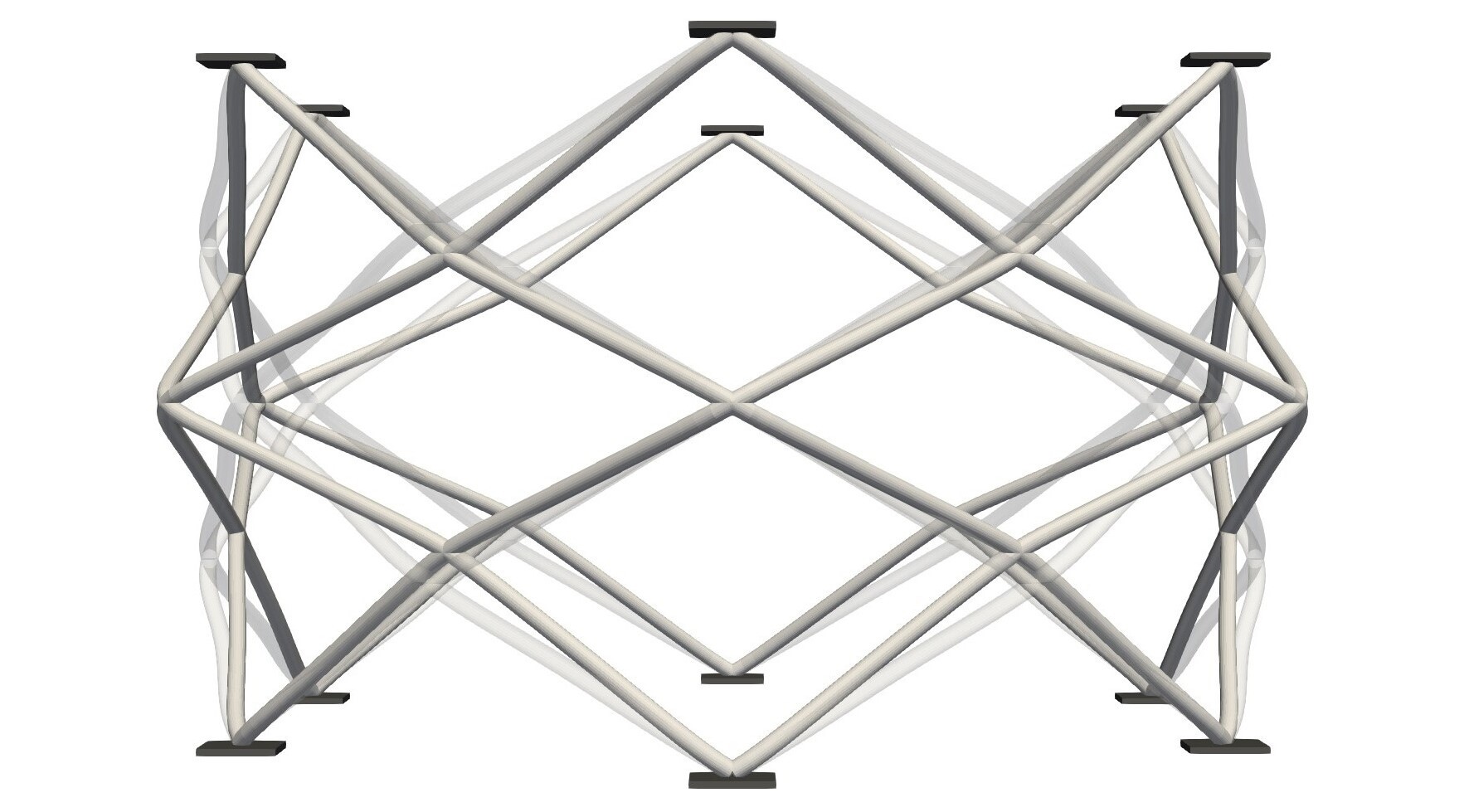}
\end{overpic}}
\subfigure[View on $x_2$-$x_3$, $t=\SI{3.00}{s}$.\label{fig:17b}]{\begin{overpic}[width=0.29\textwidth]{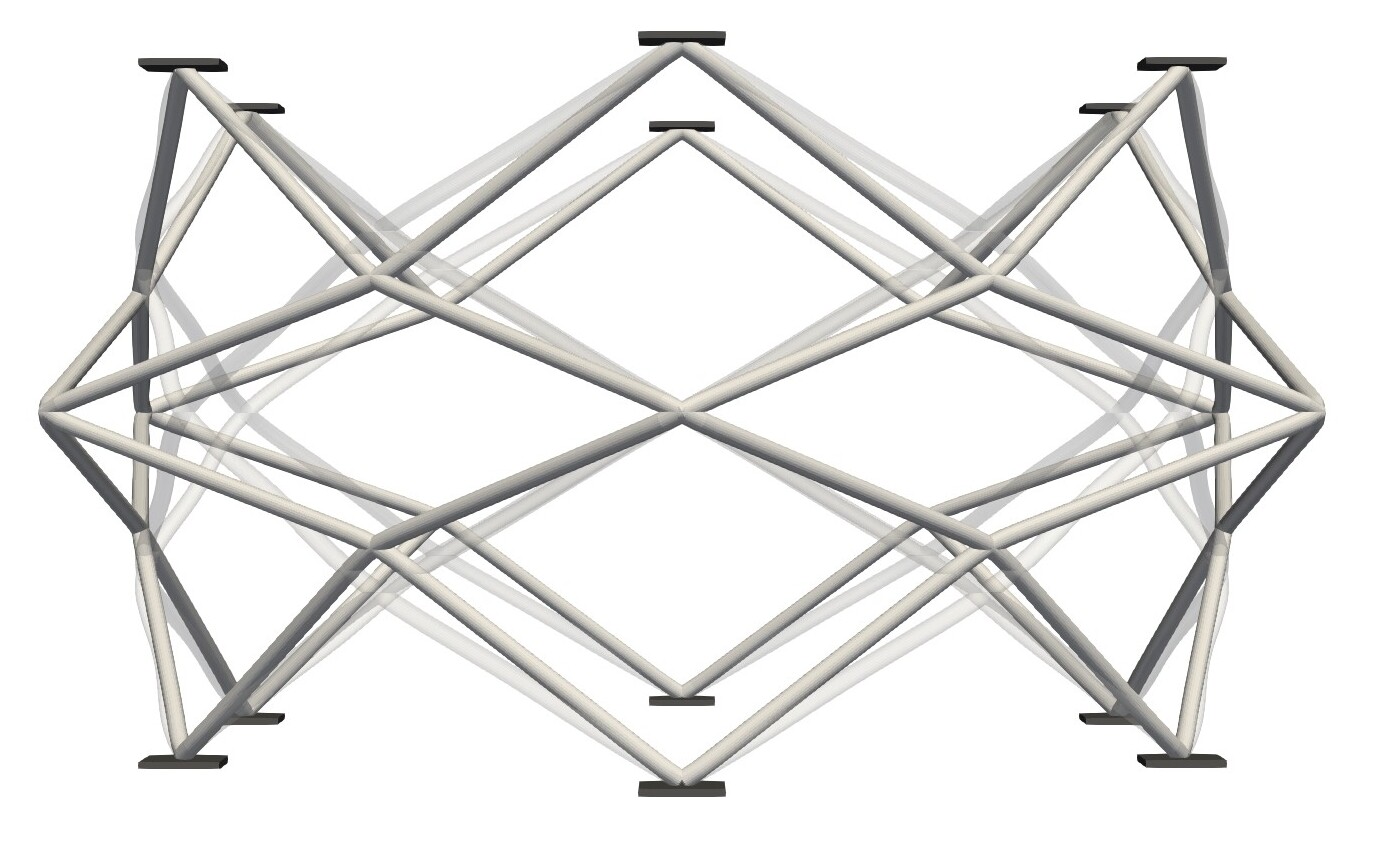}
\end{overpic}}
\subfigure[View on $x_2$-$x_3$, $t=\SI{6.00}{s}$.\label{fig:17c}]{\begin{overpic}[width=0.32\textwidth]{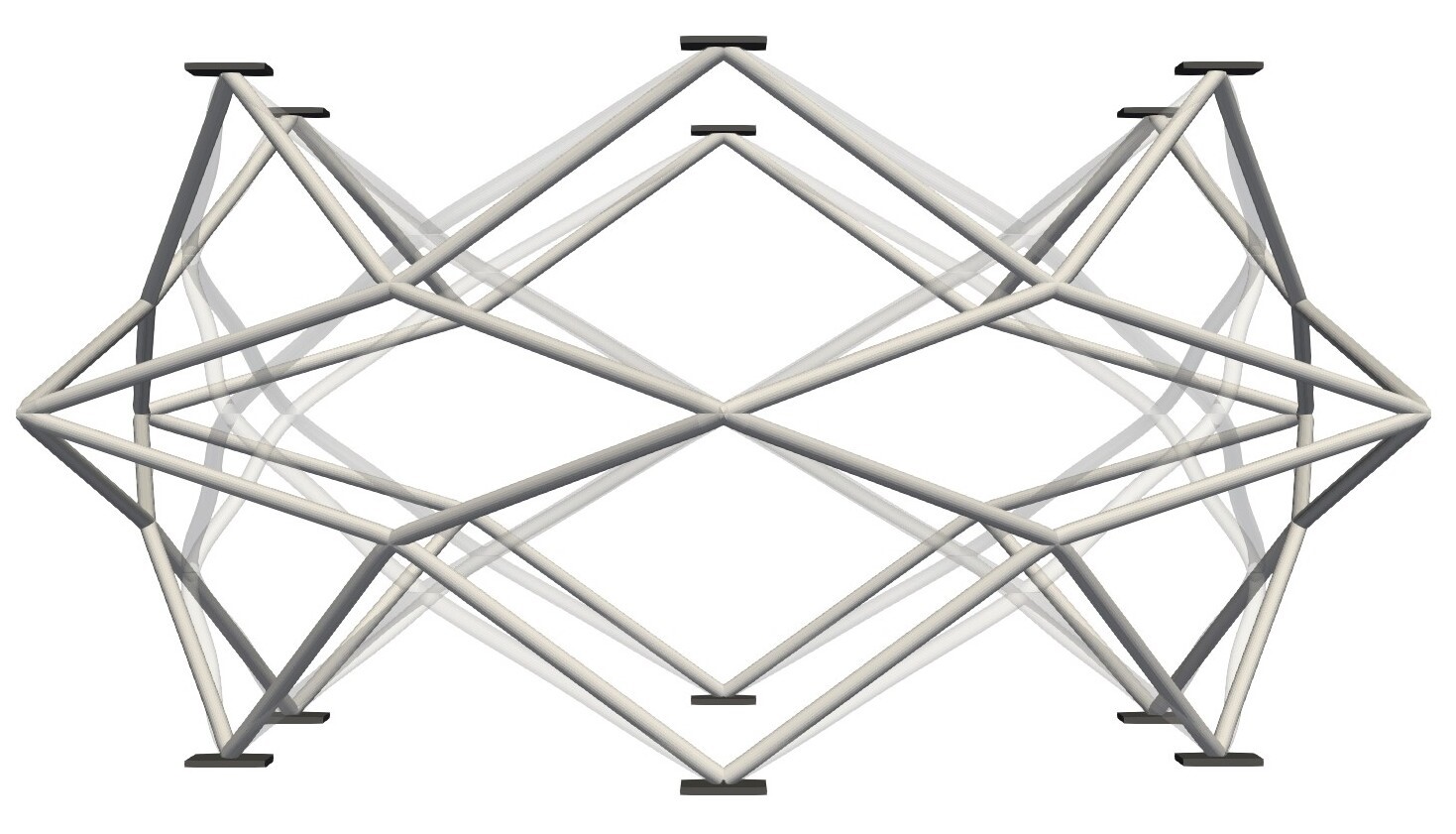}
\end{overpic}}

\subfigure[View on $x_1$-$x_3$, $t=\SI{1.15}{s}$.\label{fig:17d}]{\begin{overpic}[width=0.32\textwidth]{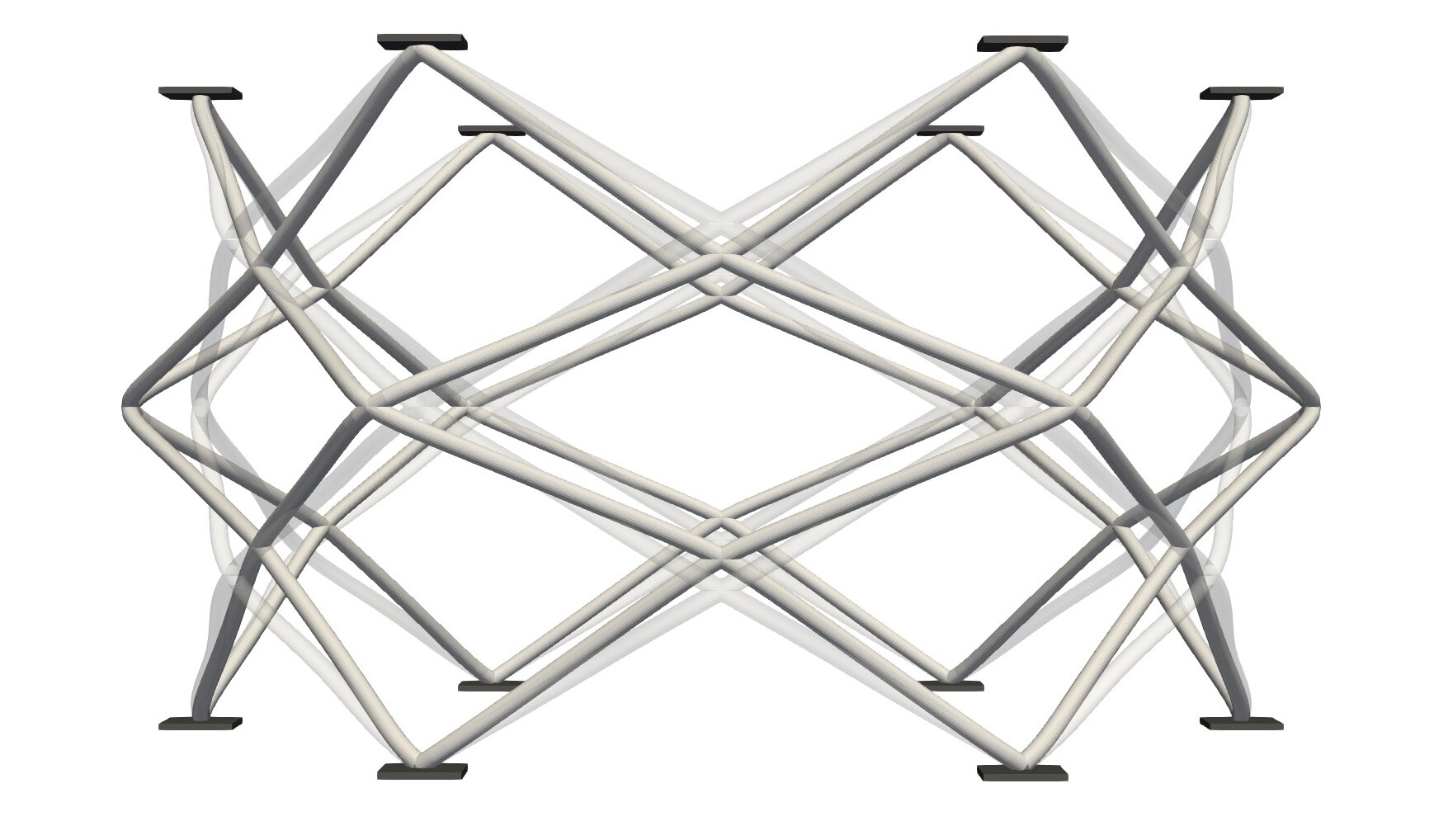}
\end{overpic}}
\subfigure[View on $x_1$-$x_3$, $t=\SI{3.00}{s}$.\label{fig:17e}]{\begin{overpic}[width=0.32\textwidth]{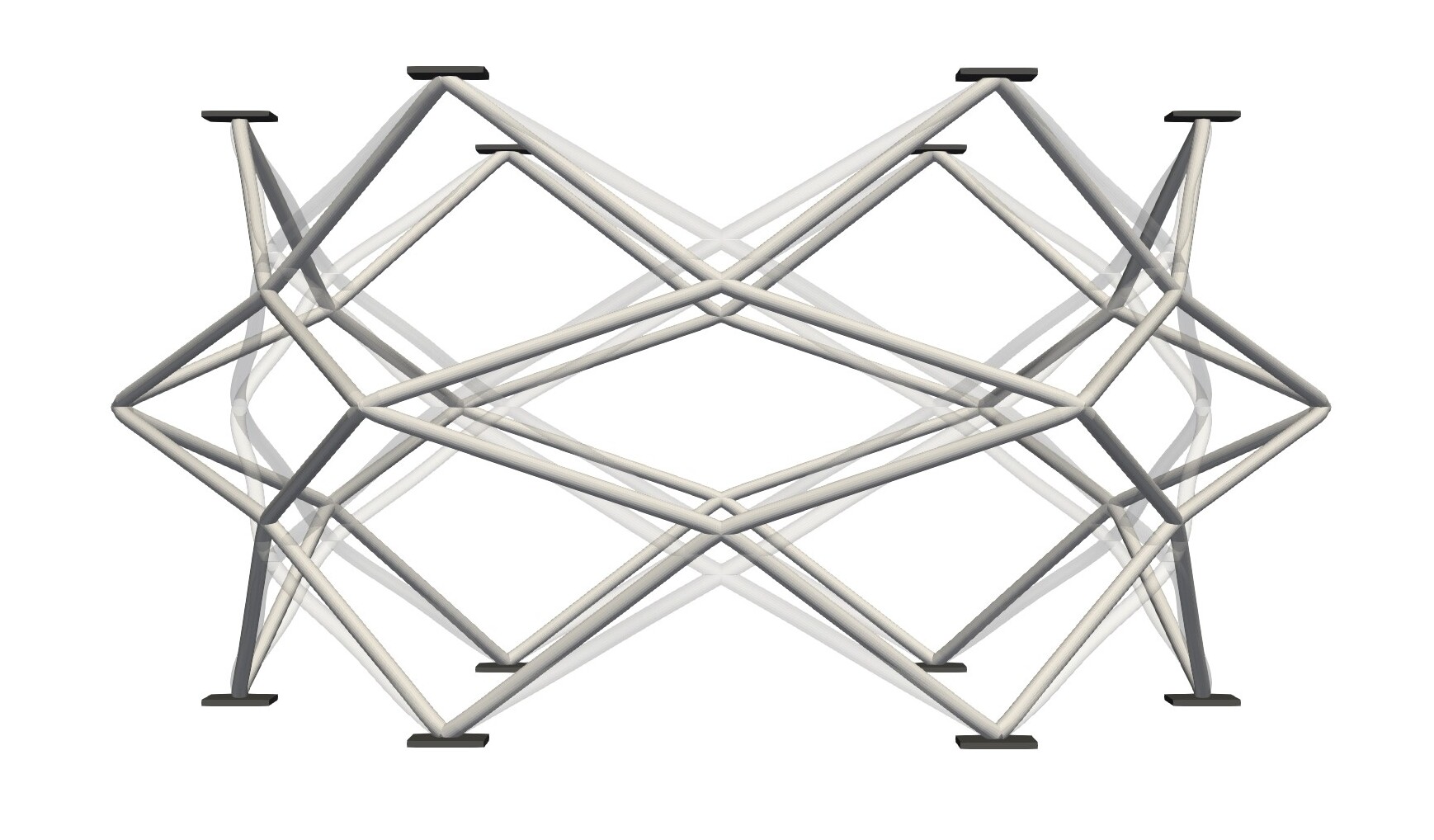}
\end{overpic}}
\subfigure[View on $x_1$-$x_3$, $t=\SI{6.00}{s}$.\label{fig:17f}]{\begin{overpic}[width=0.32\textwidth]{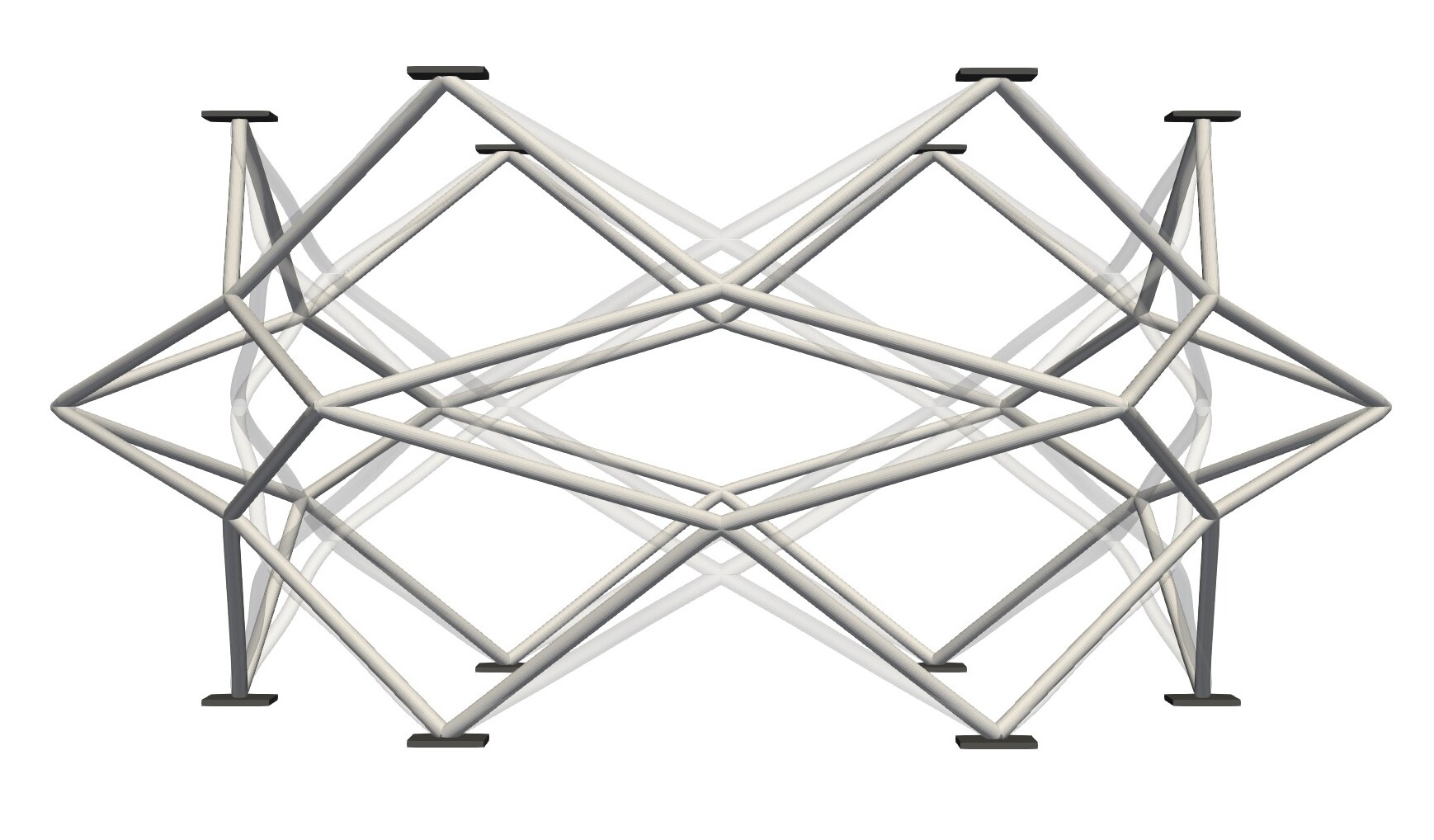}
\end{overpic}}

\subfigure[View on $x_1$-$x_2$, $t=\SI{1.15}{s}$.\label{fig:17g}]{\begin{overpic}[width=0.32\textwidth]{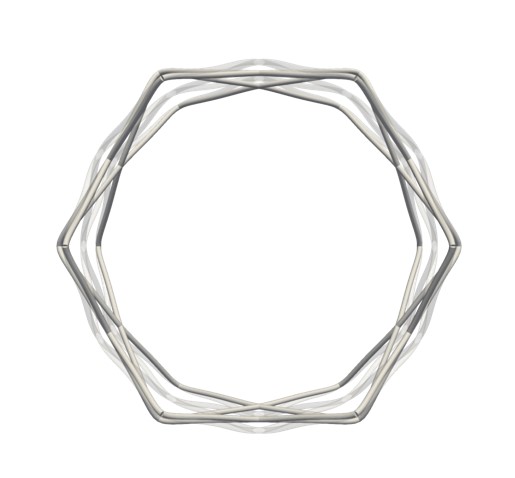}
\end{overpic}}
\subfigure[View on $x_1$-$x_2$, $t=\SI{3.00}{s}$.\label{fig:17h}]{\begin{overpic}[width=0.33\textwidth]{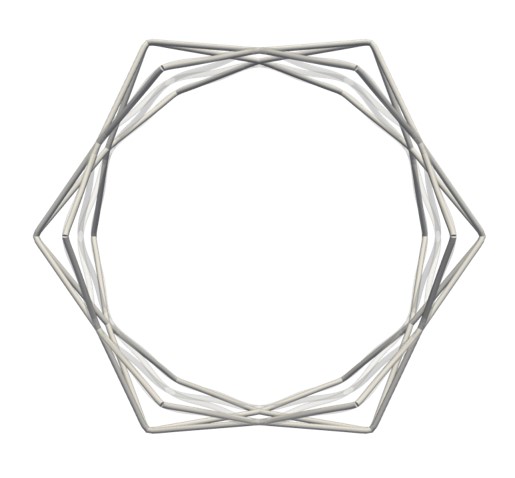}
\end{overpic}}
\subfigure[View on $x_1$-$x_2$, $t=\SI{6.00}{s}$.\label{fig:17i}]{\begin{overpic}[width=0.32\textwidth]{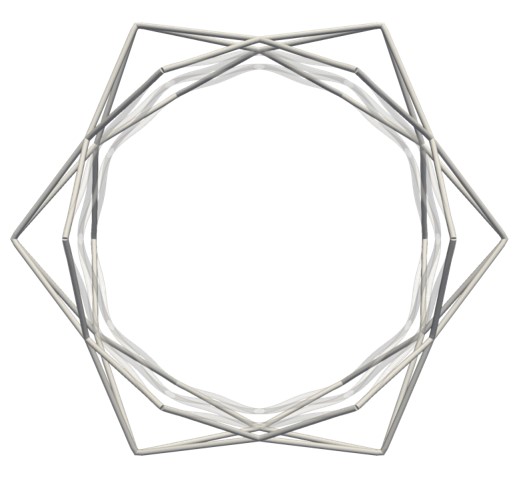}
\end{overpic}}
\caption{Tubular net: views of the deformed shape at different time instants.\label{fig:17}} 
\end{figure}

Figure~\ref{fig:17} shows snapshots of the deformation of the tubular net for $t=\SI{1.15}{s}$, \SI{3.00}{s}, and \SI{6.00}{s}. Different lateral views are shown in Figures~\ref{fig:17a}--\ref{fig:17f}, whereas top views are shown in Figures~\ref{fig:17g}--\ref{fig:17i}. The net is locally expanded due to the radial forces and increases its radius of about the 40\%. 


\section{Conclusions\label{sec:conclusions}}
Motivated by the growing demand for fast and accurate simulation tools for materials and objects with architectured inner structures, in this paper we proposed an efficient high-order formulation for geometrically exact viscoelastic beams. Linear viscoelasticity is modeled employing the generalized Maxwell model for one-dimensional solids. 
Very high efficiency is achieved by combining a number of key features. Firstly, the formulation is displacement-based, meaning that the minimal number of equations and unknowns are required. Secondly, for the spatial discretization we adopted the isogeometric collocation method, which allows to bypass integration over the elements and requires only one point evaluation per unknown. Thirdly, finite rotations are updated using the incremental rotation vector, leading to two main benefits: minimum number of rotation unknowns (namely the three components of the incremental rotation vector) and no singularity problems. Moreover, the formulation has two remarkable advantages: the same $\SO3$-consistent linearization of the governing equations and update procedures as for static or dynamic with linear elastic rate-independent materials can be directly used, and the standard second-order accurate trapezoidal rule for time integration turned out to be consistent with the underlying geometric structure of the kinematic problem. 
High-order space accuracy is obtained exploiting the IGA attributes, especially the tunable smoothness of the basis functions, the ability to accurately reconstruct complex initial geometries, and the $k$-refinement. 
Through a number of numerical applications, we demonstrated all the expected features, in particular the high-order accuracy and the robustness in managing arbitrarily complex three-dimensional beam or beams system. In our opinion, the present work opens interesting perspectives towards more efficient simulations of programmable objects, for example with applications to patient-tailored biomedical devices such as cardiovascular stents.  Next developments may include, among others, modeling thermo-responsive materials and removing the kinematic assumption of rigid beam cross sections.

\section*{Acknowledgments}
GF and EM were partially supported  by the European Union - Next Generation EU, in the context of The National Recovery and Resilience Plan, Investment 1.5 Ecosystems of Innovation, Project Tuscany Health Ecosystem (THE). (CUP: B83C22003920001).

EM was also partially supported by the National Centre for HPC, Big Data and Quantum Computing funded by the European Union within the Next Generation EU recovery plan. (CUP B83C22002830001). 

These supports are gratefully acknowledged.

\clearpage
\bibliographystyle{elsarticle-num}
\bibliography{mylibrary_new,mylibrary2_new}
\end{document}